\documentclass[11pt]{amsart}

\usepackage{enumerate}
\usepackage{latexsym}
\usepackage[pdftex]{graphicx}
\usepackage{amssymb}
\usepackage{color}

\newtheorem{theorem}{Theorem}
\newtheorem{lemma}{Lemma}
\newtheorem{proposition}{Proposition}
\newtheorem{remark}{Remark}
\newtheorem{example}{Example}

\newtheorem{definition}{Definition}
\newtheorem{corollary}{Corollary}

\newtheorem{problem}{Problem}

\def\R{{\mathbb R}}

\newcommand{\beq}{\begin{equation}}
\newcommand{\eeq}{\end{equation}}
\newcommand{\beqna}{\begin{eqnarray*}}
\newcommand{\eeqna}{\end{eqnarray*}}
\newcommand{\beqn}{\begin{equation*}}
\newcommand{\eeqn}{\end{equation*}}
\newcommand{\bp}{\begin{proof}}
\newcommand{\ep}{\end{proof}}
\newcommand{\bprop}{\begin{proposition}}
\newcommand{\eprop}{\end{proposition}}
\newcommand{\bt}{\begin{theorem}}
\newcommand{\et}{\end{theorem}}
\newcommand{\bex}{\begin{example}}
\newcommand{\eex}{\end{example}}
\newcommand{\bc}{\begin{corollary}}
\newcommand{\ec}{\end{corollary}}
\newcommand{\bl}{\begin{lemma}}
\newcommand{\el}{\end{lemma}}
\newcommand{\bprob}{\begin{problem}}
\newcommand{\eprob}{\end{problem}}
\newcommand{\br}{\begin{remark}}
\newcommand{\er}{\end{remark}}
\newcommand{\bd}{\begin{definition}}
\newcommand{\ed}{\end{definition}}

\begin{document}

\title
[Directly congruent projections]
{On   bodies    with directly congruent   projections and sections}

\author{M. Angeles Alfonseca}
\address{Department of Mathematics, North Dakota State University\\
Fargo, ND 58108, USA} \email{maria.alfonseca@ndsu.edu}

\author{Michelle Cordier}
\address{Department of Mathematics, Kent State University,
Kent, OH 44242, USA} \email{mcordier@math.kent.edu}

\author{Dmitry Ryabogin}
\address{Department of Mathematics, Kent State University,
Kent, OH 44242, USA} \email{ryabogin@math.kent.edu}

\thanks{The first author is supported in part by U.S.~National Science Foundation Grant DMS-1100657; the  third author is supported in
part by U.S.~National Science Foundation Grants DMS-0652684 and  DMS-1101636.}


\begin{abstract}
Let $K$ and $L$  be two convex bodies in ${\mathbb R^4}$, such that their projections onto all $3$-dimensional subspaces  are directly congruent.
We prove that if the set of diameters of the bodies satisfies an additional condition and some projections do not have certain $\pi$-symmetries, then 
$K$ and $L$ coincide up to  translation and an orthogonal transformation.   
We also
show that an
  analogous   statement  holds for sections of star bodies, and prove the $n$-dimensional versions of these results. 
\end{abstract}

\maketitle

\section{Introduction}

In this paper we address the following problems (see \cite[ Problem 3.2, page 125 and Problem 7.3, page 289]{Ga}).
\bprob\label{pr1} 
Suppose that $2\le k\le n-1$ and that $K$ and $L$ are convex bodies in ${\mathbb R}^n$ such that the projection $K|H$ is congruent to $L|H$ for all $H\in {\mathcal G}(n,k)$. Is $K$ a translate of $\pm L$?
\eprob
\bprob\label{pr2}
Suppose that $2\le k\le n-1$ and that $K$ and $L$ are star bodies in ${\mathbb R}^n$ such that the section $K\cap H$ is congruent to $L\cap H$ for all $H\in {\mathcal G}(n,k)$. Is $K$ a translate of $\pm L$?
\eprob
Here we say that $K|H$, the projection of $K$ onto $H$,  is congruent to $L|H$ if there exists   an orthogonal transformation  $\varphi\in O(k,H)$ in $H$ such that $\varphi(K|H)$ is  a translate of  $L|H$; ${\mathcal G}(n,k)$ stands for the Grassmann manifold of all $k$-dimensional subspaces in ${\mathbb R^n}$.

If the corresponding  projections  are translates of each other, or if the bodies are convex and the corresponding sections are translates of each other, the answers to Problems \ref{pr1} and \ref{pr2} are known to be affirmative \cite[Theorems  3.1.3 and 7.1.1]{Ga}, (see also \cite{A},  \cite{R1}). Besides, for Problem 1, with $k=n-1$, Hadwiger established a more general result and showed that it is not necessary to consider projections onto all $(n-1)$-dimensional subspaces; the hypotheses need only be true for one fixed subspace $H$, together with all subspaces containing a line orthogonal to $H$. In other words, one requires only a ``ground" projection on $H$ and all corresponding ``side" projections. Moreover, Hadwiger  noted  that in ${\mathbb R^n}$, $n\ge 4$, the ground projection might be dispensed with
(see \cite{Ha}, and \cite[pages 126--127]{Ga}).

If   the corresponding  projections (sections) of convex (star-shaped) bodies are rotations of each other,
the results in the case $k=2$ were obtained by the third author in \cite{R}; see also \cite{M}.

In the general case of rigid motions, Problems  \ref{pr1} and \ref{pr2} are open  for any $k$ and $n$. 
In the special case of {\it direct rigid motions}, {\it i.e.},  when the general orthogonal group  $O(k)$ is replaced by the special orthogonal group $SO(k)$, the problems are open as well. 

Golubyatnikov \cite{ Go1} obtained several interesting results related to Problem \ref{pr1} in the cases $k=2, 3$  \cite[Theorem 2.1.1, page 13; Theorem 3.2.1, page 48]{Go1}. 
In particular,  he  gave an affirmative answer to Problem \ref{pr1} in the  case  $k=2$ if the projections of $K$ and $L$ are directly congruent and have no direct rigid motion symmetries.

If the bodies are symmetric, then the answers to Problems 1 and 2 are known to be affirmative. In the case of projections they are
the consequence of the Aleksandrov   Uniqueness Theorem about  convex bodies, having equal volumes of projections (see \cite[Theorem 3.3.1, page 111]{Ga}); in the case of sections  they follow from
 the Generalized Funk Theorem \cite[Theorem 7.2.6, page 281]{Ga}.

In this paper we follow  the ideas from \cite{Go1} and \cite{R} to 
 obtain several Hadwiger-type results related to  both Problems 1 and 2 in the case $k=3$.
In order to formulate these results we  introduce some notation and definitions.

Let $n\ge 4$ and  let $S^{n-1}$ be the unit sphere in ${\mathbb R^n}$. We will use the notation  $w^{\perp}$ for  the 
 $(n-1)$-dimensional subspace of ${\mathbb R^n}$ orthogonal to $w\in S^{n-1}$. 
We  denote by
$d_K(\zeta)$  the diameter of  a convex body $K$, which is parallel to the direction $\zeta\in S^{n-1}$. 
We will also denote by ${\mathcal O}={\mathcal O}_{\zeta}\in O(n)$  the orthogonal transformation satisfying ${\mathcal O}|_{\zeta^{\perp}}=-I|_{\zeta^{\perp}}$, and 
${\mathcal O}(\zeta)=\zeta$. 

We define the notion of  {\it rigid motion symmetry} for sets,  
as it will be used throughout the paper. Let $D$ be a  subset of $H\in{\mathcal G}(n,k)$, $3\le k\le n-1$. We say that  $D$     has a  rigid motion symmetry if 
$\varphi(D)=D+a$
for some vector  $a\in H$ and some non-identical orthogonal transformation $\varphi\in O(k, H)$ in $H$. Similarly,   $D$  has  a direct rigid motion symmetry if $\varphi(D)=D+a$
for some vector  $a\in H$ and some non-trivial rotation $\varphi\in SO(k, H)$. 
 In the case when $D$ is a  subset of $H\in{\mathcal G}(n,3)$, and  $\xi\in (H\cap S^{n-1})$,
we say that $D$ has a 
$ (\xi,\alpha \pi)$-symmetry if $\varphi(D)=D+a$
for  some vector $a\in H$ and some rotation $\varphi\in SO(3, H)$ by the angle $\alpha\pi$, $\alpha\in (0,2)$, satisfying $\varphi(\xi)=\xi$.
If, in particular,  the angle of rotation is $\pi$, we say that $D$ has a $(\xi,\pi)$-symmetry.

\subsection{Results about directly congruent projections}
\begin{figure}[ht]
\includegraphics[clip, width=7cm]{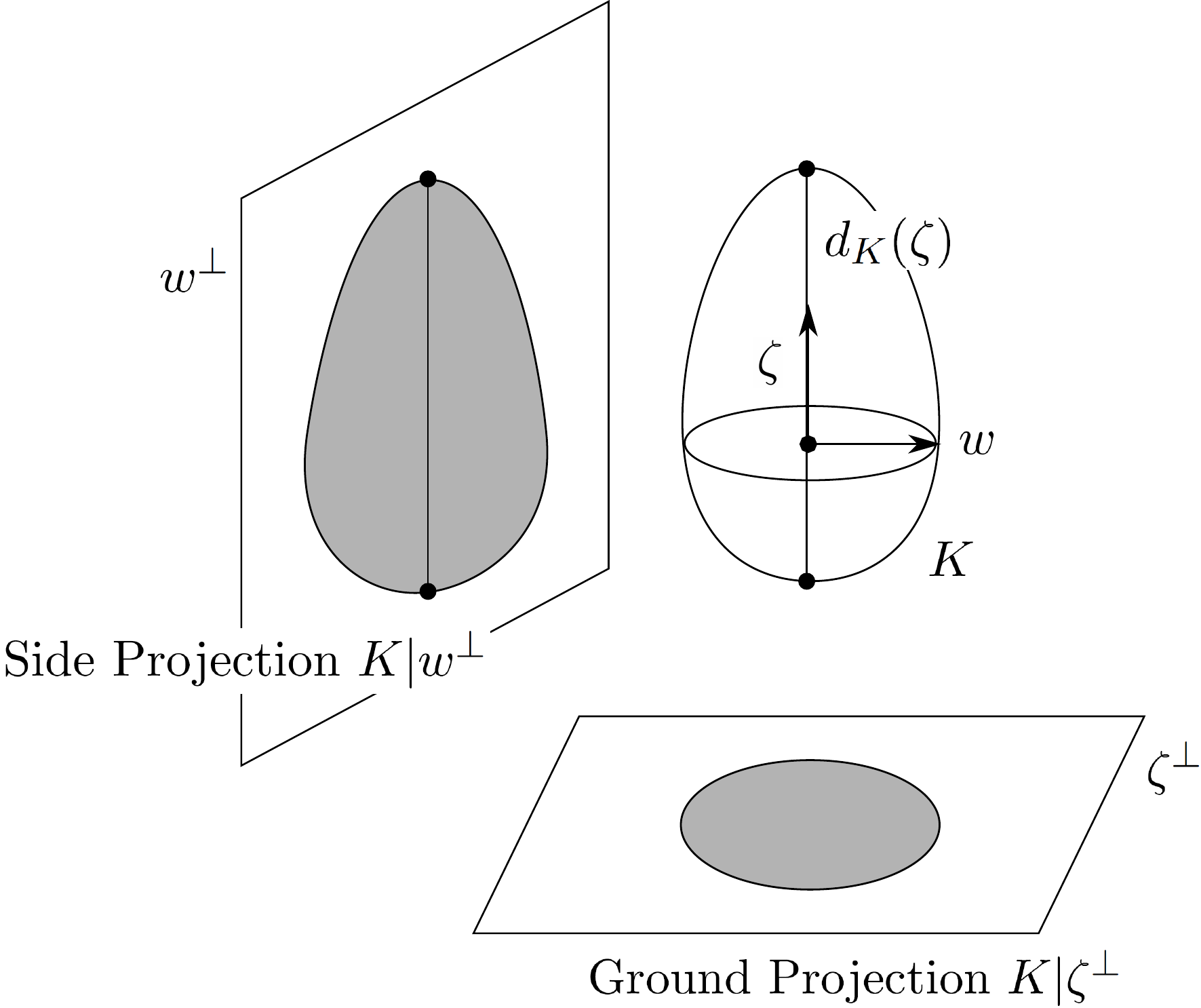}\\
\caption{Diameter $d_K(\zeta)$, side projection $K|w^\perp$ and ground projection $K|\zeta^\perp$.} 
\label{sgproj}
\end{figure}

We start with the following $4$-dimensional result.
\bt\label{tpr}
Let  $K$ and $L$ be two convex bodies in ${\mathbb R}^4$ having countably many diameters. Assume  that there exists a diameter $d_K(\zeta)$, such that the  ``side" projections $K|w^{\perp}$, $L|w^{\perp}$ onto   all subspaces
$ w^{\perp}$ containing $\zeta$   are directly congruent, see Figure \ref{sgproj}.
Assume also that these projections  have 
no $(\zeta,\pi)$-symmetries  and no $(u,\pi)$-symmetries for any $u\in (\zeta^{\perp}\cap w^{\perp}\cap  S^3)$.
Then   $K=L+b$ or $K={\mathcal O}L+b$ for some $b\in {\mathbb R^4}$.

If, in addition,  the ``ground" projections $K|\zeta^{\perp}$, $L|\zeta^{\perp}$, are directly congruent and do not have
  rigid motion symmetries,
then $K=L+b$ for some $b\in {\mathbb R^4}$.
\et

We state a straight $n$-dimensional generalization of Theorem \ref{tpr} as  a corollary.
\bc\label{npr1}
Let  $K$ and $L$ be two convex bodies in ${\mathbb R}^n$, $n\ge 4$,  having countably many diameters. Assume  that there exists a diameter $d_K(\zeta)$ such that 
the ``side" projections
$K|H$, $L|H$ onto all $3$-dimensional subspaces $H$ containing $\zeta$
  are  directly congruent. Assume also that these projections  have 
no $(\zeta,\pi)$-symmetries and no $(u,\pi)$-symmetries for any $u\in (\zeta^{\perp}\cap H\cap S^{n-1})$.
Then $K=L+b$ or $K={\mathcal O}L+b$ for some $b\in {\mathbb R^n}$.

If, in addition,   the ``ground" projections $K|G$, $L|G$ onto  all $3$-dimensional subspaces $G$ of $\zeta^{\perp}$, are  directly congruent and  have no rigid motion symmetries, then
 $ K=L+b$ for some $b\in {\mathbb R^n}$.
\ec

In particular, we see that   if 
 $K$ and $L$ are  convex  bodies in ${\mathbb R^n}$, $n\ge 4$, having countably many diameters, and directly congruent  projections onto {\bf all} $3$-dimensional subspaces, and if the ``side" and ``ground" projections related to one of the diameters satisfy the conditions of the above corollary, then $K$ and $L$ are translates of each other.

This statement was proved by Golubyatnikov  \cite[Theorem 3.2.1, page 48]{Go1} under the stronger assumptions that the ``side" projections have no direct rigid motion symmetries.
 Theorem \ref{tpr} and Corollary \ref{npr1}   under the same stronger assumptions are  implicitly contained in his proof. 
To weaken the symmetry conditions on the ``side" projections   we  replace the topological argument from \cite{Go1} with an analytic one
based on ideas from \cite{R} (compare \cite[pages 48--52]{Go1} with Proposition \ref{krh2} in Section \ref{sec2}).

We  note that the assumption  about countability of the sets of the diameters of $K$ and $L$ can be weakened.  Instead, one can assume, for example, that these sets are  subsets of a  countable union of the great circles containing $\zeta$  (see Remark \ref{ruki88} after Lemma \ref{nunu2}). We  also  note that the set of bodies considered in the above statements contains the set of all polytopes whose three dimensional projections do not have rigid motion symmetries. This set of polytopes is  an  everywhere dense set with respect to the Hausdorff metric in the class of all convex  bodies in ${\mathbb R^n}$, $n\ge 4$. For the convenience of the reader we prove  this in the Appendix.

\subsection{Results about directly congruent sections}
The analytic approach also allows  to obtain results related to Problem \ref{pr2} (see  \cite[pages 288-290, open problems 7.1, 7.3, and Note 7.1]{Ga}).
\bt\label{ts}
Let  $K$ and $L$ be two star-shaped bodies with respect to the origin in ${\mathbb R}^4$, 
having   countably many diameters. 
 Assume  that there exists  a diameter $d_K(\zeta)$ containing the origin and parallel to $\zeta$, such that for all subspaces $ w^{\perp}$ containing $\zeta$,
the  ``side" sections $K\cap w^{\perp}$, $L\cap w^{\perp}$,  are directly congruent. Assume also that these sections  have 
no $(\zeta,\pi)$-symmetries  and no $(u,\pi)$-symmetries  for any $u\in (\zeta^{\perp}\cap w^{\perp}\cap S^3)$.
Then    $K=L+b$ or $K={\mathcal O}L+b$ for some $b\in {\mathbb R^4}$ parallel to $\zeta$.

\et

As in the case of projections, we state a straight $n$-dimensional generalization of Theorem \ref{ts} as  a corollary.
\bc\label{nts1}
Let  $K$ and $L$ be star-shaped bodies with respect to the origin in ${\mathbb R}^n$, $n\ge 4$,  having countably many diameters. Assume  that there exists a diameter $d_K(\zeta)$ containing the origin, such that for all $3$-dimensional subspaces $H$ containing  $\zeta$, the ``side" sections
$K\cap H$, $L\cap H$  are  directly congruent.
Assume also that these sections  have 
no $(\zeta,\pi)$-symmetries and no $(u,\pi)$-symmetries for any $u\in (\zeta^{\perp}\cap H\cap S^3)$.
Then $ K=L+b$ or $K={\mathcal O}L+b$ for some $b\in {\mathbb R^n}$ parallel to $\zeta$.
\ec




Applying the ideas  used in this paper, one can obtain similar results  related to both Problems \ref{pr1} and \ref{pr2} in the case $k=2$, see \cite{AC}. However, we are unaware
of  results related to  the case $k\ge 4$.

The paper is organized as follows. In Section \ref{sec2} we formulate and prove our main auxiliary result,  Proposition \ref{krh2}.  Section \ref{sec3} is devoted to the  proof of  Theorem \ref{tpr} and Corollary \ref{npr1}. Theorem \ref{ts} and  Corollary  \ref{nts1} 
are proved in Section \ref{sec5}. 
We prove  that the set of polytopes in ${\mathbb R^n}$, $n\ge 4$, with $3$-dimensional projections having no rigid motion symmetries is dense in the Hausdorff metric in the class of all convex bodies in the Appendix.

\bigskip

{\bf Acknowledgements:}  We would like to thank Alexander Fish for his comments. We are very grateful to Vladimir Golubyatnikov for sharing his ideas during his visit to Kent State University. We are indebted to Mark Rudelson for suggesting the proof of Proposition \ref{AppProp}. We would like to thank the referee for his suggestions.

\section{ Notation and auxiliary definitions}

We will use the following standard notation. The unit  sphere in ${\mathbb R}^n$, $n\ge 2$, is   $S^{n-1}$. 
Given $w\in S^{n-1}$, the hyperplane orthogonal to $w$ and passing through the origin will be denoted by  $w^{\perp}=\{x\in {\mathbb R^n}:\,x\cdot w=0   \}$. Here $x\cdot w=x_1 w_1+\dots+x_n w_n$ is the usual inner product in ${\mathbb R^n}$. The Grassmann manifold of all $k$-dimensional subspaces in ${\mathbb R^n}$ will be denoted by ${\mathcal G}(n,k)$. The notation $O(k)$ and $SO(k)$, $2\le k\le n$, for the subgroups of the orthogonal group $O(n)$ and the special orthogonal group $SO(n)$  in ${\mathbb R}^n$ is standard. If ${\mathcal U}\in O(n)$ is an orthogonal matrix, we will write ${\mathcal U}^t$ for its transpose.

We refer to \cite[Chapter 1]{Ga} for the next definitions involving convex and star bodies. A {\em body} in $\R^n$ is a compact set which is equal to the closure of its non-empty interior.  A {\it convex body} is a body $K$ such that for every pair of points in $K$, the segment joining them is contained in $K$.  For $x \in \mathbb{R}^n$, the {\it support function} of a convex body $K$ is defined as 
$
h_K(x)=\max \{x\cdot y:\, \,y\in K   \}
$
 (see page 16 in \cite{Ga}). The {\it width function} $\omega_K(x)$ of $K$ in the direction $x\in S^{n-1}$ is defined as
	$\omega_K(x)=h_K(x)+h_K(-x)$. A segment $[z,y]\subset K$  is called a {\it diameter} of the convex body $K$  if $|z-y|=\max\limits_{\{\theta\in S^{n-1}\}}\omega_K(\theta)$.  We say that a convex body $K\subset {\mathbb R^n}$ has countably many diameters if the width function $\omega_K$ 
reaches its maximum on a countable subset of $S^{n-1}$.

Observe that a convex body $K$ has at most one diameter parallel to a given direction $\zeta \in S^{n-1}$ (for, if  $K$ had two parallel diameters $d_1$, $d_2$, then $K$ would contain a parallelogram with sides $d_1$ and $d_2$, one of whose  diagonals is longer than $d_1$). For this reason, if $K$ has a diameter parallel to $\zeta\in S^{n-1}$, we will denote it by  $d_K(\zeta)$.

A set $S\subset\R^n$ is said to be {\it star-shaped with respect to a point $p$}  if the line segment from $p$ to any point in $S$ is contained in $S$.  
For $x\in {\mathbb R}^n\setminus \{0\}$, and $K\subset {\mathbb R}^n$ a nonempty, compact, star-shaped set with respect to the origin, the {\it radial function} of $K$  is defined as 
$
\rho_K(x)=\max \{c:\,cx\in K   \}.
$
Here, the line through $x$ and the origin is assumed to meet $K$
(\cite[page 18]{Ga}). We say that a body $K$ is a {\it star body} if it $K$ is star-shaped with respect to the origin and its radial function $\rho_K$ is continuous.

  Given a star body $K$, a segment $[z,y]\subset K$  is called a {\it diameter} of $K$  if $|z-y|=\max\limits_{\{[a,b]\subset K\}}|a-b|$.  If a star body $K$, which is not convex, has a diameter containing the origin, that is parallel to $\zeta\in S^{n-1}$, we will also denote it by $d_K(\zeta)$.

\begin{figure}[ht]
\includegraphics[clip, width=7cm]{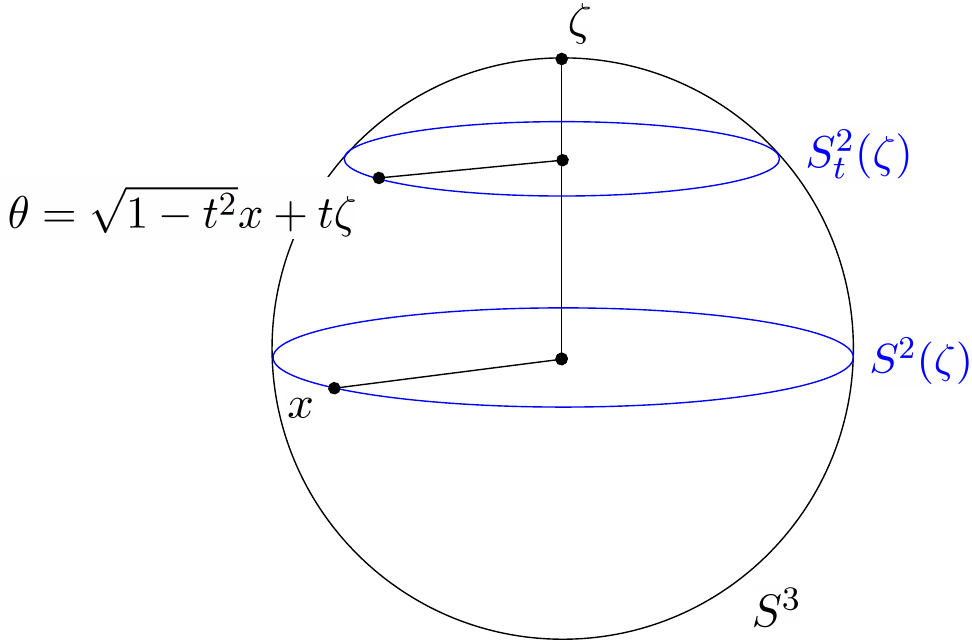}\\
\caption{The great subsphere $S^2(\zeta)$ and the parallel $S^2_t(\zeta)$.} 
\label{parallels}
\end{figure}

 Given $\zeta\in S^{n-1}$,   the great $(n-2)$-dimensional sub-sphere of $S^{n-1}$ that is perpendicular to $\zeta$ will be denoted by $S^{n-2}(\zeta)=\{\theta\in S^{n-1}:\,\theta\cdot\zeta=0   \}$.   
For $t\in [-1,1]$, the parallel to $S^{n-2}(\zeta)$ at height $t$ will  be denoted by  $S^{n-2}_t(\zeta)=S^{n-1}\cap \{x\in {\mathbb R^n}:\,x\cdot \zeta=t\}$. Observe that when $t=0$, $S^{n-2}_0(\zeta)=S^{n-2}(\zeta)$. Figure \ref{parallels} shows the case $n=4$.

 Let $E$ be a two or three-dimensional subspace  of ${\mathbb R}^n$.  We will  write $\varphi_{E}\in SO(2, E)$, or $\varphi_{E}\in SO(3, E)$, meaning that  there exists a  choice of an orthonormal basis in ${\mathbb R}^n$ and a rotation $\Phi\in SO(n)$, with a matrix written in this basis, such that  the action of $\Phi$ on
 $E$ is  the  rotation $\varphi_E$ in $E$, and the  action of $\Phi$ on  $E^{\perp}$ is trivial, {\it i.e.}, $\Phi(y)=y$ for every $y\in E^{\perp}$ (here $E^{\perp}$ stands for the orthogonal complement of $E$). A similar notation will be used for $\varphi_{E}\in O(3, E)$. For $w  \in S^3$, we will denote by $O(3, S^2(w))$, $SO(3, S^2(w))$, the orthogonal transformations in the $3$-dimensional subspace spanned by the great subsphere $S^2(w)$ of $S^3$. The restriction of a transformation $\varphi \in O(n)$ onto the subspace of smallest dimension containing  $W\subset S^{n-1}$ will be denoted by  $\varphi|_W$. $I$ stands for the identity transformation.

Finally, we define the notion of rotational symmetry for functions, as it will be used throughout the paper. Let $\zeta 
\in S^{3}$, and let  $w \in S^2(\zeta)$. For $\alpha \in [0,2]$, we will denote by $\varphi_w^{\alpha\pi}$ the rotation of the sphere $S^2(w)$ by the angle $\alpha\pi$ around $\zeta$,
{\it i.e.},  $\varphi_w^{\alpha\pi}(\zeta)=\zeta$. By this we mean that $\varphi_w^{\alpha\pi}$ is the restriction to the  3-dimensional subspace spanned by $S^2(w)$ of a  rotation $\Phi \in SO(4)$  with the following properties:  $\Phi(\zeta)=\zeta$, $\Phi(w)=w$, if  $\{x, y, w, \zeta\}$ is a positively oriented orthonormal basis of  $\mathbb{R}^4$, then for every $v \in (span\{x, y\}\cap S^3)=S^2(w)\cap S^2(\zeta)$, the angle between the vectors $v$ and  $\varphi_w^{\alpha\pi}(v)  \in S^2(w)\cap S^2(\zeta)$  is $\alpha \pi$, and if $\alpha\neq 0,1,2$, $\{v,\varphi_w^{\alpha\pi}(v),w,\zeta\}$ form a positively oriented basis of $\mathbb{R}^4$.

\bd\label{d2}
Let $f$ be a continuous function on $S^3$ and let $\xi\in S^3$. We 
 say that the restriction of $f$ onto $S^2(\xi)$ (or just $f$) has a  $(\zeta,\alpha \pi)$-rotational symmetry  if for some rotation $\varphi_{\zeta}^{\alpha\pi}\in SO(3, S^2(\xi))$ by the angle $\alpha\pi$ about  the vector $\zeta\in S^2(\xi)$, one has $f\circ \varphi_{\zeta}^{\alpha\pi}=f$ on $S^2(\xi)$. In particular, if $\alpha=1$, we say that  $f$  has a $(\zeta, \pi)$-rotational symmetry on $S^2(\xi)$.

\ed

\begin{figure}[ht]
\includegraphics[clip, width=6.5cm]{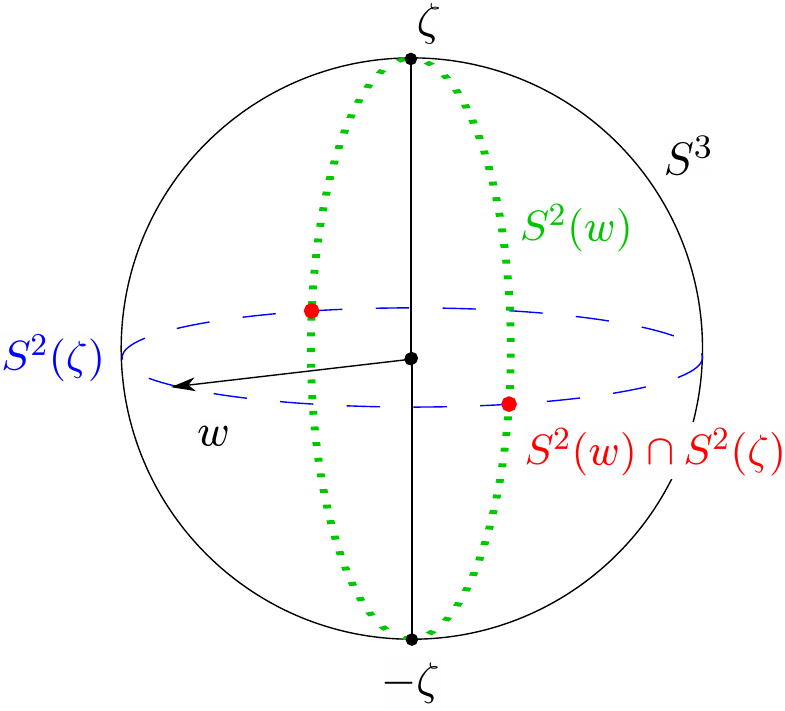}
\includegraphics[clip, width=10cm]{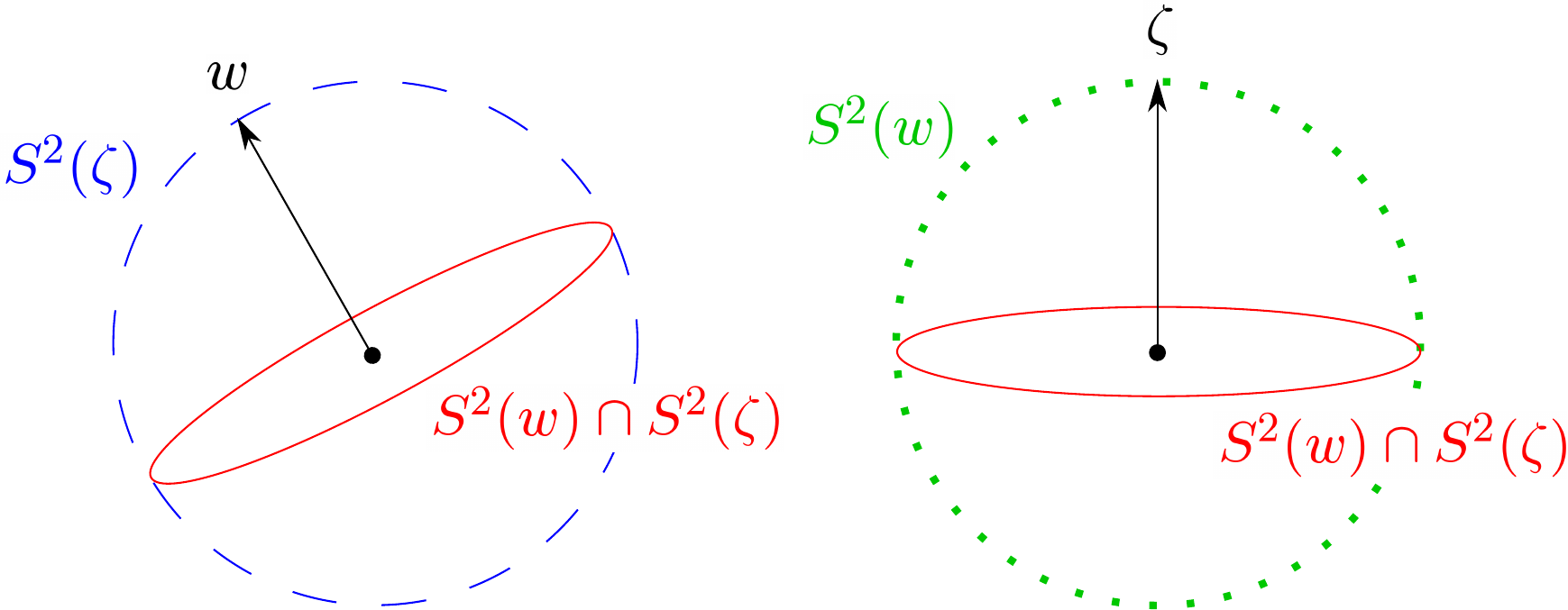}\\
\caption{The upper figure depicts $S^3$ with the two dimensional great subspheres $S^2(\zeta)$, dashed, and $S^2(w)$, dotted; the two larger dots stand for the one-dimensional subsphere $S^2(w)\cap S^2(\zeta)$. The lower figures depict $S^2(\zeta)$ and $S^2(w)$, within their corresponding 3-dimensional subspaces.}
\label{s2w}
\end{figure}

\section{A result about a  functional equation on $S^3$}\label{sec2}


In \cite{R}, the third author proved that if two continuous functions $F$ and $G$ on $S^2$ coincide up to rotation on each one-dimensional great circle, then either $F(x)=G(x)$ or $F(x)=G(-x)$ for every $x \in S^2$. The main result of this section is  a related statement for $S^3$, which, in our opinion, has independent interest.

\bprop\label{krh2}
Let $f$ and $g$ be two continuous functions on  $S^3$. Assume that for some $\zeta\in S^3$ and for every $w\in S^2(\zeta)$  there exists a rotation $\varphi_w\in SO(3, S^2(w))$, verifying that 
\begin{equation}
 \label{pos}
f\circ\varphi_w (\theta)=g(\theta), \;\;\; \forall \theta \in S^2(w).
\end{equation}
 Assume, in addition, that one of the following conditions holds:
  \begin{enumerate}[(a)]
  \item\label{pos1}   $\varphi_w(\zeta)=\zeta$   $\;\;\;\;  \forall w \in S^2(\zeta);$

  \item\label{pos1*}  $\varphi_w(\zeta)=\pm \zeta$  $\;  \forall w \in S^2(\zeta),$  
 and $f$ and $g$  have no $(\zeta,\pi)$-rotational symmetries  and no $(u,\pi)$-rotational symmetries for any $u\in S^2(\zeta)\cap S^2(w)$.
\end{enumerate}

Then either $f=g$ on $S^3$ or $f(\theta)=g({\mathcal O}\theta)$  $\forall\theta\in S^3$, where ${\mathcal O}\in O(4)$ is the orthogonal transformation satisfying ${\mathcal O}|_{S^2(\zeta)}=-I$, and 
${\mathcal O}(\zeta)=\zeta$.
\eprop

\subsection{Auxiliary Lemmata}

The direction $\zeta\in S^3$ will be fixed throughout the proof. 
We start with an easy observation about the geometry of the three dimensional sphere.


\bl\label{vse3}
Let $\zeta\in S^3$ and let $\xi\in S^2(\zeta)$. Then
\begin{equation}\label{fish1}
S^3=\bigcup\limits_{\{w\in S^2(\xi)\cap S^2(\zeta)\}}S^2(w).
\end{equation}
\el
\bp

For any $w \in S^2(\zeta)$, the two-dimensional sphere $S^2(w)$ can be written as the union of all one-dimensional parallels $S^2(w) \cap S^2_t(\zeta)$, $t \in [-1,1]$,  {\it i.e.} 
\begin{equation}\label{ruki2}
S^2(w)=\bigcup\limits_{\{t\in [-1,1] \}}(S^2(w)\cap S^2_t(\zeta)).
\end{equation}
On the other hand, we can write the two-dimensional sphere $S^2(\zeta)$ as the union of all meridians containing a fixed direction $\xi \in S^2(\zeta)$ 
as 
\[
S^2(\zeta)=\bigcup\limits_{\{w\in S^2(\xi)\cap S^2(\zeta) \}}(S^2(w)\cap S^2(\zeta)),
\]
and, rescaling, the same is true for every $S^2_t(\zeta)$, $t\in [-1,1]$ 
Thus, we have 
\begin{equation}\label{ruki1}
S^2_t(\zeta)=\bigcup\limits_{\{w\in S^2(\xi)\cap S^2(\zeta) \}}(S^2(w)\cap S^2_t(\zeta)) \qquad\forall t\in [-1,1].
\end{equation}
Combining (\ref{ruki2}) and (\ref{ruki1}), we obtain
$$
S^3=\bigcup\limits_{\{t\in [-1,1] \}} S^2_t(\zeta)=\bigcup\limits_{\{t\in [-1,1] \}}\bigcup\limits_{\{w\in S^2(\xi)\cap S^2(\zeta) \}}(S^2(w)\cap S^2_t(\zeta))=
$$
$$
\bigcup\limits_{\{w\in S^2(\xi)\cap S^2(\zeta) \}}\bigcup\limits_{\{t\in [-1,1] \}}(S^2(w)\cap S^2_t(\zeta))=
\bigcup\limits_{\{w\in S^2(\xi)\cap S^2(\zeta) \}}S^2(w).
$$
\ep


Let ${\mathcal O}\in O(4)$ be an orthogonal transformation, 
satisfying ${\mathcal O}|_{S^2(\zeta)}=-I$, and 
${\mathcal O}(\zeta)=\zeta$. Observe  that ${\mathcal O}|_{S^2(w)}$ commutes with every rotation $\varphi_w\in SO(3, S^2(w))$, such that  
$\varphi_w(\zeta)= \pm \zeta$, where $w\in S^2(\zeta)$.
It is clear that any function $f$ on $S^3$ can be decomposed in the form
\begin{equation}\label{nax2}
f(\theta)=\frac{f(\theta)+f(\mathcal O\theta)}{2}+\frac{f(\theta)-f(\mathcal O\theta)}{2}=f_{\mathcal O, e}(\theta)+f_{\mathcal O, o}(\theta),\quad\theta\in S^3,
\end{equation}
where we will call $f_{\mathcal O, e}$, $f_{\mathcal O, o}$, the even and odd parts of $f$ with respect to ${\mathcal O}$.
Since ${\mathcal O}^2=I$,  we have
$$
f_{\mathcal O, e}(\theta)=f_{\mathcal O, e}(\mathcal O\theta),\qquad f_{\mathcal O, o}(\theta)=-f_{\mathcal O, o}(\mathcal O\theta).
$$
It is also clear  that every $\theta\in S^3$ belongs to $S^2_t(\zeta)$ for some $t\in [-1,1]$, {\it i.e.},   can be written in the form
\begin{equation}\label{fru1}
\theta=\sqrt{1-t^2}x+t\zeta,
\end{equation} 
for some $t\in [-1,1]$ and $x\in S^2(\zeta)$ (see Figure \ref{parallels}).

Let $t\in [-1,1]$. For any function $f$   on $S^3$, we can define the  function $F_t$  on $S^2(\zeta)$,
\begin{equation}\label{ara}
F_t(x)=F_{t,\zeta}(x)=f(\sqrt{1-t^2}x+t\zeta),\qquad x\in S^2(\zeta),
\end{equation}
which is the restriction of $f$ to $S^2_t(\zeta)$.
Observe that 
$$
{(F_t)}_e(x)=\frac{f(\sqrt{1-t^2}x+t\zeta)+f(-\sqrt{1-t^2}x+t\zeta)}{2}=\frac{f(\theta)+f(\mathcal O\theta)}{2},
$$
where $\theta$ is as in (\ref{fru1}), {\it  i.e.},
\begin{equation}\label{uru1}
{(F_t)}_e(x)=f_{\mathcal O, e}(\theta),\qquad {(F_t)}_o(x)=f_{\mathcal O, o}(\theta).
\end{equation}
Note that  $(F_t)_e(x)=(F_t)_e(-x)$ for every $x\in S^2(\zeta)$.

As seen in the proof of Lemma \ref{vse3}, every one-dimensional great circle of $S^2(\zeta)$ is of the form $S^2(w)\cap S^2(\zeta)$ for some $w\in S^2(\zeta)$. To simplify the notation, we will denote such great circles by 
$$
E=E_{\zeta,w}=\{\theta\in S^3:\, \theta\cdot\zeta=\theta\cdot w=0  \}.
$$
Since  $\varphi_w(\zeta)=\pm \zeta$ and  $\varphi_w(S^2(w))=S^2(w)$, we have
 $$
\varphi_w(E_{\zeta,w})=\varphi_w(S^2(w)\cap S^2(\zeta))=S^2(w)\cap S^2(\zeta)=E_{\zeta,w}. 
$$


Thus, for every $t\in [-1,1]$, and for the corresponding one-dimensional equator 
$E=E_{\zeta,w}$ of $S^2(\zeta)$, there is a rotation $\phi_E \in SO(2, E)$, which is the restriction to $E$ of the  rotation  $\varphi_w\in SO(3, S^2(w))$ given by the conditions of Proposition \ref{krh2}, and which satisfies 
\begin{equation}\label{ar111}
F_t\circ\phi_E(x)=G_t(x)\qquad\forall x\in E,
\end{equation}
Here  $G_t$ is defined from $g$ similarly to $F_t$ in (\ref{ara}).

\bigskip

In the next lemma, we will need to  use the {\it Funk transform}, \cite[Chapter III, \S 1)]{He},
$$
Rf(w)=R_{\zeta}f(w)=\int\limits_{S^2(w)\cap S^2(\zeta)}f(\theta)d\theta,\qquad w\in S^2(\zeta).
$$
 Here $d\theta$ stands for the Lebesgue measure on the one-dimensional great circle $E=S^2(w)\cap S^2(\zeta)$ of $S^2(\zeta)$.

\bl\label{vse}
Let $f$ and $g$ be as in Proposition \ref{krh2}. Then $f_{\mathcal O, e}=g_{\mathcal O, e}$.
\el
\bp
Let $w\in S^2(\zeta)$, and let $\varphi_w\in SO(3, S^2(w))$ be such that (\ref{pos}) holds.
Then,
$\phi_{E}=\varphi_w|_{S^2(w)\cap S^2(\zeta)}\in SO(2, E)$ is the corresponding rotation  in  $E=S^2(w)\cap S^2(\zeta)$.
By
 the rotation invariance of the Lebesgue measure on $E$ and (\ref{ar111}), we have
 \begin{equation}\label{vasya1}
 \int\limits_{E}F_t(x)dx=\int\limits_{E}F_t\circ\phi_E(x)dx=\int\limits_{E}G_t(x)dx,\qquad \forall t\in [-1,1].
 \end{equation}
 Hence,
$R_{\zeta}F_t(w)=R_{\zeta}G_t(w)$ for every $w\in  S^2(\zeta)$.
 Thus, ${(F_t)}_e(x)={(G_t)}_e(x)$ for every $x\in S^2(\zeta)$ (apply  Theorem C.2.4 from \cite[page 430]{Ga} to ${(F_t)}_e-{(G_t)}_e$). Using the first relation in (\ref{uru1}), its analogue for $g$, and (\ref{ruki2}), we obtain the desired
 result.
\ep
\br\label{lihu1}
By the previous Lemma, to prove Proposition \ref{krh2}  we can (and from now on will) assume that the functions $f$ and $g$ are odd with respect to ${\mathcal O}$. In order to simplify the notation, from now on we will write $f$ and $g$ instead of $f_{\mathcal O, o}$ and $g_{\mathcal O, o}$. We will also write $F_t$ for $(F_t)_o$ and $G_t$ 
for $(G_t)_o$. 
\er

We consider the set $ \Psi=\{w\in S^2(\zeta):  \varphi_w(\zeta)=-\zeta\}$, and for 
any $\alpha\in [0,2]$,  the set $\Xi_{\alpha}$, defined as
\begin{equation}\label{vasya2}
\Xi_{\alpha}=\left\{w\in S^2(\zeta):\,\exists \varphi_w^{\alpha\pi}\in SO(3, S^2(w))\,\, \textrm{such that} \right.
\end{equation}
\[
\left. f\circ\varphi_w^{\alpha\pi}=g \,\,\textrm{on}\,\,  S^2(w) \textrm{ and } \varphi_w^{\alpha\pi}(\zeta)=\zeta \ \right\}.
\]

Observe that $\Xi_0=\{w\in S^2(\zeta) $ such that $f=g $ on $S^2(w) \}$, and 
\begin{equation}\label{nax1}
\Xi_1=\{w\in S^2(\zeta):\, f(\theta)=g({\mathcal O}\theta) \quad\forall\theta\in S^2(w) \}.
\end{equation}

With this notation, the result of Proposition \ref{krh2}  is that  $S^2(\zeta)=\Xi_0$ or $S^2(\zeta)=\Xi_1$, under either hypothesis \eqref{pos1} or \eqref{pos1*}.  We will divide the proof into several Lemmata.


\bl\label{vse1}
 The set  $\Xi_{\alpha}$ is closed.
\el
\bp
We can assume that ${\Xi}_{\alpha}$ is not empty.

Let $(w_l)_{l=1}^{\infty}$ be a sequence of elements of  ${\Xi}_{\alpha}$ converging to $w\in S^{2}(\zeta)$ as $l\to\infty$, and let
$\theta$ be any point on $S^2(w)$.  Consider a sequence $(\theta_l)_{l=1}^{\infty}$ of points $\theta_l\in S^2(w_l)$ converging to $\theta$ as $l\to\infty$.

(It is readily seen that such  a sequence exists.
Indeed, let $B_{\frac{1}{l}}(\theta)$ be a Euclidean ball centered at $\theta$ of radius $\frac{1}{l}$, where $l\in {\mathbb N}$.  Since $S^2(w_m)\to S^2(w)$ as $m\to\infty$, for each $l\in {\mathbb N}$ there exists $m=m(l)$ such that
$$
S^2(w_{m(l)})\cap B_{\frac{1}{l}}(\theta)\neq \emptyset.
$$
Choose any $\theta_l=\theta_{m(l)}\in S^2(w_{m(l)})\cap B_{\frac{1}{l}}(\theta)$. Then $\theta_l\to \theta$ as $l\to\infty$.)

By the definition of $\Xi_{\alpha}$, we see that
\begin{equation}\label{restl}
f\circ\varphi^{\alpha \pi}_{w_l}(\theta_l)=g(\theta_l)\qquad \theta_l\in S^2(w_l), \quad l\in {\mathbb N}.
\end{equation}
Passing to a subsequence if necessary, we can assume that the sequence of rotations $(\varphi_{w_l}^{\alpha\pi})_{l=1}^{\infty}$, 
$\varphi^{\alpha\pi}_{w_l}\in SO(3, S^2(w_l))$, is convergent, say, to $\varphi_w\in SO(3, S^2(w))$. 
Writing out the matrices of rotations $\varphi_{w_l}^{\alpha\pi}$ in the corresponding orthonormal bases $\{x_l,y_l, w_l,\zeta\}$, $x_l,y_l$$\in S^2(w_l)\cap S^2(\zeta)$, and 
passing to the limit as $l\to\infty$ we see that
$\varphi_w$ is the rotation by the angle $\alpha\pi$ and  the limit of (\ref{restl}) is
$f\circ\varphi_{w}(\theta)=g(\theta)$.
Since the choice of $\theta\in S^2(w)$ was arbitrary, we obtain that $w\in {\Xi}_{\alpha}$, and the result follows.
\ep

\br
A similar argument can be used to show that the set $\Psi$ is closed. 
\er

\bl\label{vse2}
If $\alpha\in \left({\mathbb R}\setminus{\mathbb Q}\right) \cap [0,2]$, then $\Xi_{\alpha}\subset \Xi_0$.
\el
\bp
Assume that $\Xi_\alpha \neq \emptyset$ and let $w\in \Xi_{\alpha}$. Following the ideas of Schneider \cite{Sch1}, we claim at first that $f^2=g^2$ on $S^2(w)$.
Indeed, since $f$ and $g$ are odd with respect to ${\mathcal O}$, $f^2$ and $g^2$ are even with respect to ${\mathcal O}$, and satisfy \eqref{pos} with $f^2$, $g^2$ instead of $f$, $g$. Thus, by Lemma \ref{vse}, we obtain        that  $f^2=g^2$ on $S^2(w)$. 

Squaring (\ref{pos}), we have (with $\varphi_w=\varphi_w^{\alpha \pi}$), 
$$
f^2\circ\varphi_w(\theta)=g^2(\theta)=f^2(\theta)\qquad \forall\theta\in S^2(w).
$$
Iterating for any $k\in{\mathbb Z}$,
$$
f^2\circ\varphi^k_w(\theta)=f^2\circ\varphi^{k-1}_w(\theta)=\dots=f^2(\theta)\qquad \forall \theta\in S^2(w),
$$
and using the fact that for every $\theta\in S^2(w)$, the orbit of $(\varphi^k_w(\theta))_{k\in{\mathbb Z}}$ is dense on every parallel $ S^2(w)\cap S^2_t(\zeta)$ of $S^2(w)$ 
(where $t \in [-1,1]$), we obtain that 
the restrictions of $f^2$ and $g^2$ onto $S^2(w)$ are invariant under rotations leaving $\zeta$ fixed. In other words, $f^2$ and $g^2$ are constant on every parallel of $S^2(w)$ orthogonal to $\zeta$. By continuity, $f$ and $g$ must also be constant on these parallels and $f\circ\varphi_w=f$. Hence, using (\ref{pos}) we have $f=g$ on $S^2(w)$, and therefore $w\in \Xi_0$.  Since $w$ from $\Xi_{\alpha}$ was chosen arbitrarily,  we obtain  the desired result.
\ep

In Lemma \ref{vse2}, we have shown that rotations whose angle is an irrational multiple of $\pi$ are not relevant under the assumptions of Proposition \ref{krh2}. Our next goal is to prove that rational multiples are not relevant either, except for the rotations by the angles $0$ and $\pi$. This will be achieved in Lemma \ref{vse4}, by means of a topological argument, which is based on one definition and two Lemmata from \cite{R} (see Lemmata \ref{mainl} and \ref{musya1} below).  The argument will show that for each $t \in (-1,1)$ and an appropriate $w\in S^2(\zeta)$, the subset of a great circle $S^2(w)\cap S^2(\zeta)$, where the functions $F_t=G_t$ are equal to each other, is open. Since such a set is closed by definition, and it is non-empty, we will conclude that $F_t$ equals $G_t$ on this large circle.  Using  (\ref{ruki2}) we will obtain that $f=g$  on  the corresponding $S^2(w)$, which will give us the desired result.

\begin{figure}[ht]
\includegraphics[clip, width=10cm]{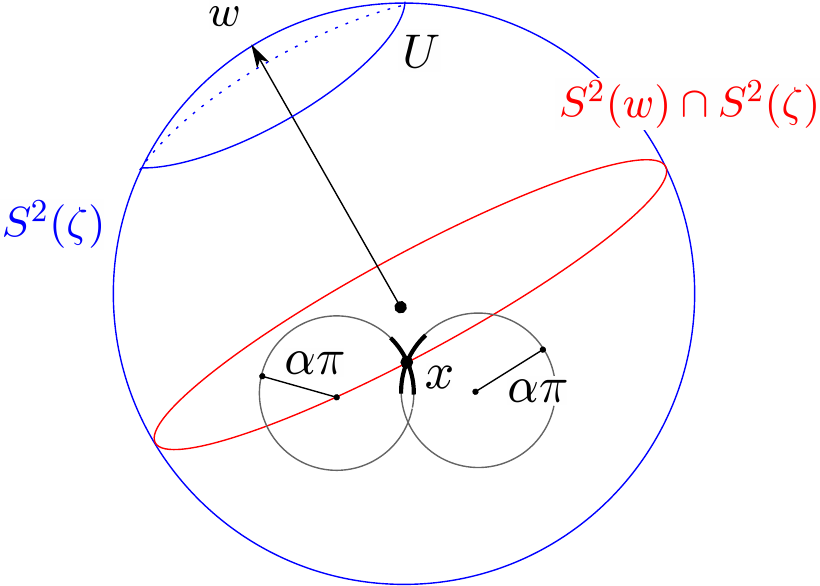}\\
\caption{The spherical $X$-figures from Definition \ref{ruba1}}
\label{X-figure}
\end{figure}

We will  reformulate the corresponding statements from \cite{R} in a way that is  more convenient for us here. Refer to Figure \ref{X-figure} for the next definition.

\bd\label{ruba1}
 Let $\alpha\in (0,1)$ and let ${\mathbf S}_1$, ${\mathbf  S}_2$ be any two spherical circles in the standard metric of $S^2(\zeta)$,  both of radius $\alpha\pi$. The union   ${\mathfrak l}\cup  {\mathfrak m}$ of two {\it open} arcs ${\mathfrak l}\subset {\mathbf S}_1$ and $ {\mathfrak m}\subset {\mathbf S}_2$ will be called a {\it spherical   $X$-figure} if  the angle between the arcs is in $(0,\frac{\pi}{4})$, the length of the arcs is less than $\alpha\pi$,  and the arcs intersect at their centers only,  ${\mathfrak l}\cap  {\mathfrak m}=\{x\}$.
The point  $x\in S^2(\zeta)$ will be called the {\it center} of the $X$-figure.
\ed

Let  $t\in (-1,1)$,   $F_t$ be a function on $S^2(\zeta)$,  and  $x$ be the center of a spherical $X$-figure. If for every $u\in X$ we have $F_t(u)=F_t(x)$, we will say that there exists an $X$-figure $X_{F_t(x)}\subset S^2(\zeta)$.
The following two Lemmata are Lemma 10 and Lemma 12  from \cite[pages 3438--39]{R} (with $f=F_t$, $g=G_t$, $f_e=F_t^2$, $S^2=S^2(\zeta)$ and $S^2(w)\cap S^2(\zeta)$ instead of $\xi^{\perp}$).

\bl\label{mainl}
Let $\,\,t\in (-1,1)$, and let $F_t$ and $G_t$  be two continuous functions on $S^2(\zeta)$.  Assume that there is an open  spherical cap $U\subset \Xi_{\frac{p}{q}}$, with $\frac{p}{q}\in (0,1)\cap{\mathbb Q}$,  such that for every  $ w\in U$, there exists a rotation $\phi_w=\phi_{w,\zeta}$ of the great circle $S^2(w)\cap S^2(\zeta)$ by the angle $\frac{p}{q}\pi$, verifying  
\begin{equation}\label{ar1}
F_t\circ\phi_{w}(x)=G_t(x) \qquad\forall x\in S^2(w)\cap S^2(\zeta).
\end{equation}
Then, for every $ x\in S^2(w)\cap S^2(\zeta)$  there exists an $X$-figure  $ X_{F_t^2(x)}\subset S^2(\zeta)$,
with  one of the arcs of $X_{F_t^2(x)}$ being orthogonal to $S^2(w)\cap S^2(\zeta)$.
Moreover, for every $ x,y \in  S^2(w)\cap S^2(\zeta)$ there exist $X$-figures $ X_{F_t^2(x)},\,X_{F_t^2(y)}\in S^2(\zeta)$, such that  
\[
\Theta(X_{F_t^2(x)})=X_{F_t^2(y)},
\]
where $\Theta\in SO(3, S^2(\zeta))$ is such that $\Theta(w)=w$ and
$\Theta(x)=y$.
\el

\bl\label{musya1}
Let $t\in (-1,1)$, and let $F_t$, $G_t$, and $U$  be as above. Then,  for every $w\in U$  there exists a constant $c$ such that
$F_t^2(x)=G_t^2(x)=c$ for every $x\in S^2(w)\cap S^2(\zeta)$.
\el

Observe that since any two great circles of $S^2(\zeta)$ intersect, the above constant is actually independent of $w\in U$.


\bl\label{lihu2}
Let $t\in (-1,1)$, and let $F_t$, $G_t$, and $U$  be as above. Then $f=g=0$ on $S^2(w)$  for every $ w\in U$.
\el
\bp
Let $w$ be any point in $U$, and let $ t \in (-1,1)$. By Remark \ref{lihu1}, $f$ and $g$ are odd with respect to ${\mathcal O}$ on $S^3$. Using the second relation in (\ref{uru1}), we see that $F_t$, $G_t$ are odd on $S^2(\zeta)$. By continuity, there exist $x_1$, $x_2\in  S^2(w)\cap S^2(\zeta)$ such that 
$F_t(x_1)=G_t(x_2)=0$. 
By Lemma \ref{musya1}, $F_t^2(x)=G_t^2(x)=0$ for every $x\in S^2(w)\cap S^2(\zeta)$.
Using  (\ref{ara}) and  the continuity of $f$ and $g$, we see that the last statement is true for all $t \in [-1,1]$. 
Finally, using  (\ref{ruki2}) and (\ref{ara}) again, we conclude that $f=g=0$ on $S^2(w)$. 
\ep
Now we are ready to prove
\bl\label{vse4}
We have $S^2(\zeta)=\Xi_0\cup\Xi_1 \cup \Psi$.
\el
\bp
Assume that the set $A:=S^2(\zeta)\setminus(\Xi_0\cup\Xi_1 \cup \Psi)$ is not empty.  By  Lemma \ref{vse2}, $A\cap \Xi_{\alpha}=\emptyset$, provided that $\alpha\in [0,2] \setminus{\mathbb Q}$. Hence, $A$ may be written as
 $$
A=\bigcup\limits_{\frac{p}{q}\in{\mathbb Q}\cap [0,2]} \left( A \cap \Xi_{\frac{p}{q}}\right).
$$
By Lemma \ref{vse1}, all $\Xi_{\frac{p}{q}}$ are closed and $A$ is open. Hence, by the Baire category Theorem, (cf. Lemma 8 from \cite{R}),
 there exists
$\frac{p}{q}\in {\mathbb Q \cap [0,2]}$ such that $int(\Xi_{\frac{p}{q}})\neq\emptyset$.
We can assume that there exists an open spherical cap $U\subseteq(A\cap \Xi_{\frac{p}{q}})$ such that for every  $w\in U$, there is a rotation $ \varphi_w^{\frac{p}{q}\pi}\in SO(3, S^2(w))$ such that 
$$
  f\circ\varphi_w^{\frac{p}{q}\pi}=g \,\,\textrm{on}\,\,  S^2(w).
$$
In particular,  for any $t\in (-1,1)$, and for every large circle $E=S^2(w)\cap S^2(\zeta)$ of $S^2(\zeta)$ there exists a rotation $\phi_w\in SO(2,  E)$ by the angle $\frac{p}{q}\pi$ such that (\ref{ar1})  holds. Changing the orientation if necessary, we can assume that $p/q$ is between 0 and 1. 

By Lemma \ref{musya1}, $F_t^2(x)=G_t^2(x)=c$ for every $x\in S^2(w)\cap S^2(\zeta)$, and by Lemma \ref{lihu2} we have
$f=g= 0$ on $S^2(w)$. Hence, $w\in \Xi_0$, which is impossible, since $w\in A$.
The result follows.
\ep

\subsection{Proof of Proposition \ref{krh2} under hypothesis (\ref{pos1})} Under hypothesis \eqref{pos1}, the set $\Psi$ is empty. Thus, by Lemma \ref{vse4}, we have that $S^2(\zeta)=\Xi_0 \cup \Xi_1$.

The proof of the next result is similar to the one of Lemma 1 in \cite[page 3434]{R} (with $\Xi_1$ instead of $\Xi_\pi$ and $S^2(\xi)\cap S^2(\zeta)$ instead of $\xi^{\perp}$).

\bl\label{ruki3}
Let $\zeta\in S^3$, $\xi\in S^2(\zeta)$. Assume that 
\[
(S^2(\xi)\cap S^2(\zeta))\cap  \Xi_0\cap \Xi_1=\emptyset. 
\]
Then, either 
 \begin{equation}\label{br5}
(S^2(\xi)\cap S^2(\zeta))\subset (\Xi_0\setminus\Xi_1)\quad \textrm{or}\quad (S^2(\xi)\cap S^2(\zeta))\subset (\Xi_1\setminus\Xi_0).
\end{equation}
\el

\bigskip

If we assume that  $\Xi_1=\emptyset$, then  $S^2(\zeta)=\Xi_0$, and therefore   $f(\theta)=g(\theta)$  for every $\theta\in S^3$. On the other hand,  if   $\Xi_0=\emptyset$, we have that $S^2(\zeta)=\Xi_1$, which means that  
$f(\theta)=g({\mathcal O}\theta)$ for every $\theta\in S^3$. Hence, in these two situations we obtain the desired conclusion. 

Let us now assume that both $\Xi_0$, $\Xi_1$ are not empty.  We can also assume that $\Xi_0\cap \Xi_1\neq\emptyset$. Indeed, let $w$ be a point on the boundary of $\Xi_0$, ($w\in\Xi_0$, since $\Xi_0$ is closed). Then for every $l\in{\mathbb N}$, the set $B_{\frac{1}{l}}(w)\cap S^3$ contains a point $w_l$ from $\Xi_1$. But then $w_l\to w$ as $l\to\infty$, hence $w\in\Xi_1$, and $w\in \Xi_0\cap \Xi_1$.

We shall consider two cases:

1) There exists  $\xi\in S^2(\zeta)$ such that $\Xi_0\cap\Xi_1\cap  S^2(\xi)= \emptyset$.

2)  For every  $x\in S^2(\zeta)$ we have $\Xi_0\cap\Xi_1\cap S^2(x)\neq \emptyset$.

Consider the first case.  Using Lemma \ref{ruki3}, we obtain (\ref{br5}).
 If the first relation in (\ref{br5}) holds, then, by Lemma \ref{vse3}, we have $S^3=\bigcup\limits_{\{w\in\Xi_0\}}S^2(w)$, and
 $f(\theta)=g(\theta)$ for every $\theta\in S^3$. 
 If the second relation in (\ref{br5}) holds, then, using Lemma \ref{vse3} again, we obtain $S^3=\bigcup\limits_{\{w\in\Xi_1\}}S^2(w)$, and
$f(\theta)=g({\mathcal O}\theta)$ for every $\theta\in S^3$.

Consider the second case. We claim  that
\begin{equation}\label{pusto}
S^2(\zeta)= \bigcup\limits_{\{u\in (\Xi_0\cap \Xi_1)\}} (S^2(u)\cap S^2(\zeta)).
\end{equation}
Indeed, let $x\in S^2(\zeta)$. By the hypothesis of the second case, the set $\Xi_0\cap\Xi_1\cap S^2(x)$ is non-empty. Let $u \in \Xi_0\cap\Xi_1\cap S^2(x)$. Then  $x \in S^2(u)$, and  hence $x \in S^2(u)\cap S^2(\zeta)$, from which it follows that 
\[
   x\in  \bigcup\limits_{\{u\in (\Xi_0\cap \Xi_1)\}} (S^2(u)\cap S^2(\zeta)),  
\]
thus proving (\ref{pusto}).
Using (\ref{pusto}), the fact that $S^3=\bigcup\limits_{\{v\in S^2(\zeta)\}}S^2(v)$, and an argument similar to the one in  the proof of Lemma \ref{vse3},  we conclude that 
\begin{equation}\label{pusto1}
S^3= \bigcup\limits_{\{u\in (\Xi_0\cap \Xi_1)\}} S^2(u).
\end{equation}


It is easy to see that if (\ref{pusto1}) holds, then
$f$ and
$g$ are {\it zero} on $S^3$, and we are done. 
Indeed, let $\theta\in S^3$. Then $\theta\in S^2(w)$ for some $w\in (\Xi_0\cap \Xi_1)$. Using (\ref{nax1}) we see that
$f(\theta)=g(\theta)=g({\mathcal O}\theta)$.
Since $g$ is odd with respect to ${\mathcal O}$, we have $g(\theta)=f(\theta)=0$. Since $\theta$ was arbitrary, we 
 have proved that if (\ref{pusto1}) holds, then $f=g=0$ on $S^3$. 

Thus, in all possible cases, we have shown that if  $f$ and $g$ are odd with respect to ${\mathcal O}$, then either
$f(\theta)=g(\theta)$ for every $\theta\in S^3$, or
$f(\theta)=g({\mathcal O}\theta)$ for every $\theta\in S^3$  (see Remark \ref{lihu1}). Proposition \ref{krh2} is proved under hypothesis \eqref{pos1}.
\qed




\subsection{Proof of Proposition \ref{krh2} under hypothesis (\ref{pos1*})}

By Lemma \ref{vse4}, we have that $S^2(\zeta)=\Xi_0 \cup \Xi_1\cup \Psi$.  We will show that the additional hypothesis on the lack of symmetries for $f$ and $g$ implies that $\Psi$ must be empty.  
This will be achieved in Lemmata \ref{m1111}$-$\ref{m5}.

\bl\label{m1111} 
We have $(\Xi_0\cup\Xi_1)\cap \Psi= \emptyset$.
\el
\bp



Assume that $(\Xi_0\cup\Xi_1)\cap \Psi$ is nonempty, and let $w  \in (\Xi_0\cup\Xi_1)\cap \Psi$. Using  the definition of $\Xi_0, \Xi_1$ and $\Psi$, we have
 $$
f\circ\varphi_w=g,\qquad f\circ\psi_w=g\qquad \textrm{on}\quad S^2(w),
$$
 where  $\varphi_w$, $\psi_w \in SO(3, S^2(w))$ 
are rotations satisfying
$\varphi_w(\zeta)=\zeta$,  $\psi_w(\zeta)=-\zeta$.

 If $w\in \Xi_0$, then
$\varphi_w$ is trivial, and we have $f= f\circ\psi_w$ on $S^2(w)$. 
Since any $3$-dimensional rotation has a one-dimensional invariant subspace, there exists $u\in S^2(w)\cap S^2(\zeta)$ such that $\psi_w(u)=u$. This means that $f$ has a $(u,\pi)$-rotational symmetry,  which is impossible by the  assumptions of Proposition  \ref{krh2}.

If $w\in \Xi_1$, then $\varphi_w$ is the rotation by angle $\pi$ around $\zeta$, while  $\psi_w$ is the rotation by  angle $\pi$ around $u\in S^2(w)\cap S^2(\zeta)$. Since $\varphi_w^{-1}=\varphi_w$, it follows that  $f= f\circ \varphi_w\circ\psi_w$.  It is well known  (see, for example, \cite{RS}) that the composition of two rotations by $\pi$ about axes that are separated by an angle $\beta$,  is a rotation by $2\beta$ about an axis perpendicular to the axes of the given rotations. Since $\zeta$ and $u$ are perpendicular,  we conclude that $\varphi_w\circ\psi_w$ is a rotation by $\pi$ around $v\in S^2(w)\cap S^2(u)\cap S^2(\zeta)$. Hence, $f$ has a $(v,\pi)$-rotational symmetry, which is impossible by the assumptions of Proposition \ref{krh2}.  Thus, $(\Xi_0\cup\Xi_1)\cap \Psi=\emptyset$, and the Lemma is proved.
\ep

Lemma \ref{m1111} implies  that either $S^2(\zeta)=\Xi_0 \cup \Xi_1$ or $S^2(\zeta)=\Psi$. The first case was already considered in the proof of Proposition \ref{krh2} with hypothesis \eqref{pos1}. 

Consider the second case. To prove the next lemma we will need  the following result of Radin and Sadun \cite{RS}: Let ${\mathcal A}$ and ${\mathcal B}$ be rotations of finite order in the Euclidean $3$-space, about axes that are themselves separated by an angle which is  a rational multiple of $\pi$. Then, the $2$-generator subgroup of $SO(3)$, generated by ${\mathcal A}$ and ${\mathcal B}$, is infinite and dense, except in the following cases: if one generator has order $1$, the group  is cyclic; if one generator has order $2$ and the axes are orthogonal, the group is dihedral; and if both generators have order $4$ and the axes are orthogonal, the group is the symmetries of the cube.

\bl\label{m555} 
 Assume that $S^2(\zeta)=\Psi$. Then $\forall w \in S^2(\zeta)$ there exists a unique  rotation $\varphi_w\in SO(3,S^2(w))$  by the angle $\pi$ around  some $u\in S^2(w)\cap S^2(\zeta)$, satisfying \eqref{pos}.
\el
\bp
Assume that for some $w\in S^2(\zeta)$ there exist  two different  rotations, $\tilde{\varphi_1}\neq \tilde{\varphi_2}$,  around $u_1\neq \pm u_2$, $u_1, u_2 \in S^2(w)\cap S^2(\zeta)$, satisfying
\begin{equation}\label{2xuxuxo1}
f\circ\tilde{\varphi_1}(\theta)=g(\theta),\quad 
 f\circ\tilde{\varphi_2}(\theta)=g(\theta)\qquad\forall\theta\in S^2(w).
\end{equation}
Then, $f\circ\tilde{\varphi_1}(\theta)= f\circ\tilde{\varphi_2}(\theta)$ for every $\theta\in S^2(w)$. 
 This implies that 
$f= f\circ\tilde{\varphi_1}\circ\tilde{\varphi_2}$ on $S^2(w)$,
where $\tilde{\varphi_1}\circ \tilde{\varphi_2}$ is the rotation by angle $2\beta$ around $\zeta$, and $\beta$ is the angle between $u_1$ and $u_2$. 
Hence, $f$ has a $(\zeta,2\beta)$-rotational symmetry. We claim that this is impossible.

Hypothesis \eqref{pos1*} excludes the case $\beta=\frac{\pi}{2}$.  If $\beta$  is   a rational multiple of $\pi$, $\beta\neq \frac{\pi}{2}$, by the remarks before Lemma \ref{m555}, we see that the  $2$-generator subgroup, generated  by 
 $\tilde{\varphi_1}$, $\tilde{\varphi_2}$ (both of which have order 2)  is dense in $SO(3,S^2(w))$. 
Using (\ref{2xuxuxo1}), we obtain
\begin{equation}\label{33xux}
f^2\circ\tilde{\varphi_1}(\theta)=g^2(\theta)=f^2\circ\tilde{\varphi_2}(\theta)\qquad\forall\theta\in S^2(w), 
\end{equation}
which implies
\begin{equation}\label{3xux}
f^2\circ\tilde{\varphi_1}(\theta)=f^2\circ\tilde{\varphi_2}(\theta)=f^2\circ\tilde{\varphi_1}\circ\tilde{\varphi_2}(\theta)=f^2(\theta)\qquad\forall\theta\in S^2(w).
\end{equation}
Since for every $\theta\in S^2(w)$ the sequence of points  $\tilde{\varphi_1}(\theta)$, $\tilde{\varphi_2}(\theta)$, 
$\tilde{\varphi_1}\circ\tilde{\varphi_2}(\theta)$,...,  generated by the words with letters $\tilde{\varphi_1}$, $\tilde{\varphi_2}$, is dense in $S^2(w)$, the functions
$f^2$ and $g^2$ must be identically constant (zero, since they are the squares of the odd functions) on $S^2(w)$, and hence $f$ and $g$ must equal zero. Then, $w\in \Xi_0 \cap \Psi$, which 
contradicts Lemma \ref{m1111}.

Similarly, if  $\beta$  is   an irrational multiple of $\pi$,  then, using (\ref{3xux}) and an argument similar to the one in Lemma \ref{vse2}, we see that $f$ is constant on every parallel of $S^2(w)$ orthogonal to $\zeta$.  We conclude that  $f$ has a $(\zeta,\pi)$-symmetry. 

Thus,  the rotation  $\varphi_w\in SO(3,S^2(w))$ must be unique, and the lemma is proved.
\ep

The idea of the proof of the following statement is taken from \cite[Lemma 3.2.1, page 48, and the third paragraph on page 51]{Go1}.
\bl\label{m5} 
Assume that $S^2(\zeta)=\Psi$. Then there exist a continuous tangent  line field  on $S^2(\zeta)$.
\el
\bp

Let  ${\mathbb A}$ be the function assigning to each $w \in S^2(\zeta)$ the  rotation ${\mathbb A}(w)=\varphi_w\in SO(3,S^2(w))$  by the angle $\pi$ around  some $u\in S^2(w)\cap S^2(\zeta)$,    $\varphi_w(\zeta)=-\zeta$. By Lemma \ref{m555}, the map 
${\mathbb A}$ is well-defined.   We claim that ${\mathbb A}$  is continuous. Let $(w_l)_{l=1}^{\infty}$ be a convergent sequence of directions from $S^2(\zeta)$, with $\lim\limits_{l\to\infty}w_l=w$, and let $(\varphi_l)_{l=1}^{\infty}$ be the corresponding sequence of  rotations in  $S^2(w_l)$,with  $\varphi_l(\zeta)=-\zeta$,  for every $ l\in \mathbb{N}$.
First, we prove that $(\varphi_l)_{l=1}^{\infty}$ is  convergent. Let     $(\theta_l)_{l=1}^{\infty}$, $\theta_l\in  S^2(w_l)$, be a sequence   converging to any point $\theta\in S^2(w)$ as $l\to\infty$ (the existence of such a sequence can be shown as in Lemma \ref{vse1}).
If $(\varphi_l)_{l=1}^{\infty}$ were not convergent, then there would exist two subsequences $(\varphi_{m_l})_{l=1}^{\infty}$ and  $(\varphi_{j_l})_{l=1}^{\infty}$, 
with $\tilde{\varphi_1}:=\lim\limits_{l\to\infty}\varphi_{m_l}\neq $$\lim\limits_{l\to\infty}\varphi_{j_l}=:\tilde{\varphi_2}$.
 Passing to the limit  in the  equalities
$$
f\circ\varphi_{w_{m_l}}(\theta_{m_l})=g(\theta_{m_l}),\qquad 
 f\circ\varphi_{w_{j_l}}(\theta_{j_l})=g(\theta_{j_l}),
$$
on the corresponding equators $S^2(w_{m_l})$, $S^2(w_{j_l})$, and using the fact that $\theta$ was an arbitrary point in $S^2(w)$,  we obtain (\ref{2xuxuxo1}). As we saw in the proof of the previous lemma, this is impossible. This contradiction shows that  the sequence $(\varphi_l)_{l=1}^{\infty}$ is convergent. 

To show that ${\mathbb A}$ is continuous, it remains to prove that  $\lim\limits_{l\to\infty}\varphi_l=\varphi_w$.
Assume that the last equality is not true, and let $\lim\limits_{l\to\infty}\varphi_l=\tilde{\varphi_1}\neq \varphi_w$.
Then  we have (\ref{2xuxuxo1}) with $\tilde{\varphi_2}=\varphi_w$, which is, as we have already seen,  impossible.
Thus, ${\mathbb A}$ is continuous.

Consider now the map ${\mathbb B}$ assigning to each $w\in S^2(\zeta)$ the one-dimensional invariant subspace ${\mathcal Y}(w)$ of the corresponding rotation $\varphi_w\in SO(3, S^2(w))$, $\varphi_w(\zeta)=-\zeta$.  A similar argument to  the one used for ${\mathbb A}$ allows us to show that the map ${\mathbb B}$ is well-defined and continuous. Observe also that   ${\mathcal Y}(w)\subset (w^{\perp}\cap \zeta^{\perp})$. Thus,  assuming that  $S^2(\zeta)=\Psi$, we have constructed  a {\it continuous tangent  line field} ${\mathcal Y}(w)$ on $S^2(\zeta)$. 
\ep

It is a  well-known result of Hopf (see \cite{Mi}, \cite{Sa}) that if a compact differentiable manifold $M$ admits a  continuous tangent line field, then the Euler characteristic of $M$ is zero. Since the Euler characteristic of the two-dimensional sphere is $2$, the assumption that $S^2(\zeta)=\Psi$ leads to a contradiction, as seen in Lemma \ref{m5}. Proposition \ref{krh2} is proven.

\section{Proofs of Theorem \ref{tpr} and Corollary \ref{npr1}}\label{sec3}

The proof of Theorem \ref{tpr} relies on the existence of a diameter $d_K(\zeta)$ of the body $K$, for which we have information about the side projections of $K$ and $L$. The main idea of the proof \cite{Go1} is to observe that if the width function $\omega_K$ achieves its maximum at the direction $\zeta$, then the hypotheses of Theorem \ref{tpr}  imply that the body $L$ also has a diameter in the direction $\zeta$, and both diameters have the same length (Lemma \ref{nunu2}). Therefore, we can translate the bodies to make the diameters coincide and be centered at the origin (Lemma \ref{naxuek1}).  Next, since $K$ and $L$ have countably many diameters, it follows that almost all 3-dimensional projections of the translated bodies $\tilde{K}$ and $\tilde{L}$ contain only this particular diameter, and thus the direct rigid motion given by the statement of Theorem 1 must fix it. There are only two possibilities, namely, that the rigid motion is a rotation around the diameter, or a rotation around a line perpendicular to the diameter.  We thus reduce matters to Proposition \ref{krh2} with $f=h_{\tilde{K}}$ and $g=h_{\tilde{L}}$.

\subsection{Auxiliary Lemmata}
 \label{auxl}
 
Let  $\zeta\in S^3$ be the direction of the diameter $d_K(\zeta)$ given by the statement of Theorem \ref{tpr}. 
By hypothesis,  the projections $K|w^{\perp}$ and $L|w^{\perp}$ are directly congruent for every $ w\in S^2(\zeta)$.
Hence, for every $w\in S^2(\zeta)$ there exists $\chi_w\in SO(3,S^2(w))$ and $a_w\in w^{\perp}$ such that
\begin{equation}\label{nax21}
\chi_w(K|w^{\perp})=L|w^{\perp}+a_w.
\end{equation}

Let ${\mathcal A}_{K}\subset S^3$ be the set of directions  parallel to the diameters of $K$, and ${\mathcal A}_{L}\subset S^3$ be the set of directions  parallel to the diameters of $L$. We define
\begin{equation}\label{nax3}
\Omega=\{w\in S^2(\zeta):\quad ({\mathcal A}_K\cup {\mathcal A}_L)\cap S^2(w)=\{\pm\zeta\}\}.
\end{equation}

We will repeatedly use  the following well-known properties of the support function. For every convex body $\tilde{K}$,  
\begin{equation}\label{nax35}
h_{\tilde{K}|w^{\perp}}(x)=h_{\tilde{K}}(x)   \;  \mbox{ and } \;  h_{\chi_w(\tilde{K}|w^{\perp})}(x)=h_{\tilde{K}|w^{\perp}}(\chi_w^t(x)), \quad \forall x \in w^{\perp},
\end{equation}
(see, for example, \cite[ (0.21), (0.26), pages 17--18]{Ga}).

Our first goal is  to reduce matters to rotations fixing the one-dimensional subspace containing $\zeta$.
We will do this by showing that for  most of the directions $w\in S^2(\zeta)$
the projections $K|w^{\perp}$ and $L|w^{\perp}$ have exactly one diameter, parallel to $\zeta$.

\bl\label{nunu2}
Let $K$ and $L$ be as in Theorem  \ref{tpr}, and let
$\zeta \in {\mathcal A}_K$. 
 Then    $\zeta \in {\mathcal A}_L$ and $\Omega$ is everywhere dense in $S^2(\zeta)$. Moreover, for every $w\in \Omega$ we have $\chi_w(\zeta)=\pm\zeta$ and $\omega_K(\zeta)=\omega_L(\zeta)$.
\el
\bp
Using (\ref{nax35}), we see  that the length of diameters $d_{K|w^{\perp}}(\zeta)$ and $d_K(\zeta)$ is the same for every $w\in S^2(\zeta)$.
Let $\xi$ be any element of ${\mathcal A}_L$, and let  $w\in S^2(\zeta)$ be such that $S^2(w)\ni \zeta, \xi$.
Since 
 $K|w^{\perp}$ and $L|w^{\perp}$ are directly congruent, and  the length of the diameters is not changed under rigid motions, we have that the diameter of $K$ in the direction $\zeta$ has the same length as the diameter of $L$ in the direction $\xi$.
 

We  prove  that $\Omega$ is everywhere dense in $S^2(\zeta)$.  Suppose $\xi \in (S^2(\zeta) \setminus \Omega)$. Then there exists $\eta \in \left( {\mathcal A}_K \cup {\mathcal A}_L \right) \cap S^2(\xi)$,  $\eta \not=\pm \zeta$. 
Hence, $\xi \in S^2(\eta) \cap S^2(\zeta)$ and 
$$
(S^2(\zeta)\setminus \Omega)\subseteq\bigcup\limits_{\{\eta\in {\mathcal A}_K\cup {\mathcal A}_L,\eta \not= \pm \zeta\}} \left( S^2(\eta) \cap S^2(\zeta) \right).
$$
Since  the right-hand side of the above inclusion is a countable union of one-dimensional circles,  the measure of 
$S^2(\zeta)\setminus \Omega$ is zero. 
Hence, $\Omega$ is everywhere dense in $S^2(\zeta)$.

We now show  that  $\zeta\in {\mathcal A}_K$ implies  $\zeta\in {\mathcal A}_L$. 
By definition of $\Omega$,  we have ${\mathcal A}_K\cap S^2(w)=\{\pm\zeta\}$ for every $w\in \Omega$. 
If $ {\mathcal A}_L\cap S^2(w)=\emptyset$, then the projection of $K$ on the subspace containing $S^2(\zeta) \cap S^2(w)$ contains the diameter of $K$, and the corresponding projection of $L$ does not. Therefore,  the width functions satisfy 
 $\omega_L(\theta)<\omega_K(\zeta)$ for every $\theta\in S^2(w)$. This  contradicts the fact that $K|w^{\perp}$ and $L|w^{\perp}$ are  directly congruent.
Thus, $ {\mathcal A}_L\cap S^2(w)=\{\pm\zeta\}$. 

Finally, assume that for some $w\in \Omega$ we have $\chi_w(\zeta)\neq \pm\zeta$. Then $\chi_w(K|w^{\perp})$ has a diameter in a direction $\eta\neq\pm\zeta$. Since $\chi_w(K|w^{\perp})$ and $L|w^{\perp}$ are translations of each other, $L|w^{\perp}$ must have a diameter parallel to $\eta$, which is impossible. Hence for every, $w\in \Omega$ we have $\chi_w(\zeta)= \pm\zeta$, and $\omega_K(\zeta)=\omega_L(\zeta)$.
The result follows.
\ep

\br\label{ruki88}
The previous lemma remains valid if,  instead of the condition about countability of the diameters of the bodies,
one assumes that, say,  the sets of diameters of $K$ and $L$ are  countable unions  of  large circles
 containing $\zeta$. 
The only fact that was used in the proof is   that the set  of the directions $w\in S^2(\zeta)$,  such that $d_K(\zeta)$ and $d_L(\zeta)$ are the only diameters of the  projections $K|w^{\perp}$ and $L|w^{\perp}$, is dense in $S^2(\zeta)$. 
\er

Our next goal is to ``separate" translations from rotations. We translate the bodies $K$ and $L$ by vectors $a_K$, $a_L\in{\mathbb R}^4$,  to obtain $\tilde{K}=K+a_K$ and $\tilde{L}=L+a_L$ such that their diameters $d_{\tilde{K}}(\zeta)$ and $d_{\tilde{L}}(\zeta)$ coincide and are centered at the origin.
\bl\label{naxuek1}
Let $\chi_w$ be the rotation given by (\ref{nax21}), and let $w\in \Omega$. Then the rotation $\varphi_w:=\left(\chi_w\right)^{t}$ 
satisfies $\varphi_w(\zeta)=\pm\zeta$ and 
\begin{equation}\label{naxuek2}
h_{\tilde{K}}\circ \varphi_w(\theta)=h_{\tilde{L}}(\theta)\qquad \forall\theta\in S^2(w).
\end{equation}
\el
\bp
Define
$b_w=\chi_w(a_K|w^{\perp})-a_L|w^{\perp}+a_w$, 
where $a_K|w^{\perp}$, $a_L|w^{\perp}$ are the projections of the vectors $a_K$, $a_L$, onto $w^{\perp}$.
Then  (\ref{nax21}) holds with $\tilde{K}$ and $\tilde{L}$ instead of $K$ and $L$, and $b_w$ instead of $a_w$. 
We claim at first that $b_w=0$ for all $w\in \Omega$. In other words,
\begin{equation}\label{nax33}
\chi_w(\tilde{K}|w^{\perp})=\tilde{L}|w^{\perp}.
\end{equation}
Indeed,  using the definition of $\tilde{K}$ and $\tilde{L}$, and  Lemma \ref{nunu2}, for every $w\in \Omega\subset S^2(\zeta)$ we have 
\[
d_{\tilde{K}|w^{\perp}}(\zeta)=d_{\tilde{K}}(\zeta)=d_{\tilde{L}}(\zeta)=
d_{\tilde{L}|w^{\perp}}(\zeta),\quad
  \chi_w(d_{\tilde{K}}(\zeta))=d_{\tilde{K}}(\zeta).
\]
It follows that 
$$
d_{\tilde{K}|w^{\perp}}(\zeta)=\chi_w(d_{\tilde{K}|w^{\perp}}(\zeta))=d_{\tilde{L}|w^{\perp}}(\zeta)+b_w=d_{\tilde{K}|w^{\perp}}(\zeta)+b_w.
$$
Thus, $b_w=0$ and (\ref{nax33}) holds for every $w\in \Omega$. 
  Then,
$$
h_{\chi_w(\tilde{K}|w^{\perp})}(x)=h_{\tilde{L}|w^{\perp}}(x) \qquad
\forall x\in w^{\perp},  
$$
together with (\ref{nax35})  gives us the desired conclusion. 
\ep

\subsection{Proof of Theorem \ref{tpr}} 

Consider the closed sets $\Xi=\{w\in S^2(\zeta):  (\ref{naxuek2}) \textrm{ holds with }\varphi_w(\zeta)=\zeta\}$ and $\Psi=\{w\in S^2(\zeta): (\ref{naxuek2}) \textrm{ holds with } \varphi_w(\zeta)=-\zeta\}$.
Since the set $\Omega\subset(\Xi\cup \Psi)$ is everywhere dense in $S^2(\zeta)$ by   Lemma \ref{nunu2}, we have that $\Xi\cup \Psi= S^2(\zeta)$.  We have thus reduced matters to Proposition \ref{krh2} with $f=h_{\tilde{K}}$ and $g=h_{\tilde{L}}$.  
Therefore, either $h_{\tilde{K}}=h_{\tilde{L}}$ on $S^3$  or $h_{\tilde{K}}(\theta)=h_{\tilde{L}}({\mathcal U}\theta)$ for every $\theta\in S^3$, where ${\mathcal U}\in O(4)$ is the orthogonal transformation satisfying ${\mathcal U}|_{S^2(\zeta)}=-I$, and ${\mathcal U}(\zeta)=\zeta$. Letting  ${\mathcal O}={\mathcal U}^t$, it follows from (\ref{nax35}) that  $h_{\tilde{K}}({\mathcal U}\theta)=h_{{\mathcal O}\tilde{K}}(\theta)$ for every $\theta\in S^3$, and  either  $K+a_K=L+a_L$ or  $K+a_K={\mathcal O}L+{\mathcal O}(a_L)$. This proves the  first part of the Theorem. 

Assume, in addition,  that the ground projections $K|\zeta^{\perp}$, $L|\zeta^{\perp}$, are directly congruent. Then, there exists $\chi_{\zeta}\in SO(3,S^2(\zeta))$ and $a_{\zeta}\in\zeta^{\perp}$ such that
$
\chi_{\zeta}(K|{\zeta}^{\perp})=L|{\zeta}^{\perp}+a_{\zeta}.
$
If $K={\mathcal O}L+b$ holds, then we have 
$$
K|{\zeta}^{\perp}=({\mathcal O}L)|{\zeta}^{\perp}+b|{\zeta}^{\perp}=-L|{\zeta}^{\perp}+b|{\zeta}^{\perp}.
$$
Therefore,
$\chi_{\zeta}(K|{\zeta}^{\perp})-a_{\zeta}=-K|{\zeta}^{\perp}+b|{\zeta}^{\perp}$, and  $K|{\zeta}^{\perp}$ has a rigid motion symmetry, which contradicts our assumptions. We conclude that  $K=L+b$, and the proof of  Theorem \ref{krh2} is finished.
\qed

\subsection{Proof of Corollary \ref{npr1}}

Let $J$ be an arbitrary  $4$-dimensional subspace of ${\mathbb R^n}$, containing $\zeta$.
Observe that $K|J$ and $L|J$ satisfy the conditions of Theorem \ref{tpr} with $K|J$ and $L|J$ instead of $K$ and $L$. We translate the bodies $K$ and $L$ by vectors $a_K$, $a_L\in{\mathbb R}^n$, to obtain $\tilde{K}=K+a_K$ and $\tilde{L}=L+a_L$ such that the origin is the center of $ d_{\tilde{K}}(\zeta)=d_{\tilde{L}}(\zeta)$. 

By Theorem \ref{tpr} we have $\tilde{K}|J=\tilde{L}|J$ or $\tilde{K}|J={\mathcal O}_J(\tilde{L}|J)$ where ${\mathcal O}_J \in O(4,J)$,  ${\mathcal O}_J|_{\zeta^\perp}=-I$ and ${\mathcal O}_J(\zeta)=\zeta$. 
If there existed two different $4$-dimensional subspaces  $J_1$ and $J_2$, such that $\tilde{K}|J_1=\tilde{L}|J_1$ and $\tilde{K}|J_2={\mathcal O}_{J_2}(\tilde{L}|J_2)$, then  $\tilde{L}$ would have a 3-dimensional projection with a $(\zeta,\pi)$-symmetry. Indeed, assume that $J_1 \cap J_2$ is a 3-dimensional subspace. 
Then, 
$$\tilde{L}|(J_1\cap J_2)=(\tilde{L}|J_1)|(J_1\cap J_2)=(\tilde{K}|J_1)|(J_1\cap J_2)=(\tilde{K}|J_2)|(J_1\cap J_2)$$
$$=({\mathcal O}_{J_2}(\tilde{L}|J_2))|(J_1 \cap J_2)={\mathcal O}_{J_2}|_{J_1}(\tilde{L}|(J_1\cap J_2)),$$ 
and 
 $\tilde{L}|(J_1\cap J_2)$
has a $(\zeta,\pi)$-symmetry, contradicting the assumptions of the Corollary. Hence, either  $\tilde{K}|J=\tilde{L}|J$ for every $J$, or  $\tilde{K}|J={\mathcal O}_J(\tilde{L}|J)$ for every $J$.  If we are in the second case, let ${\mathcal O} \in O(n)$ such that ${\mathcal O}|_{\zeta^\perp}=-I$ and  ${\mathcal O}(\zeta)=\zeta$.  Then we have  that ${\mathcal O}|_J={\mathcal O}_J$.  Since $J$ was arbitrary,  the projections of $\tilde{K}$ and $\tilde{L}$
onto all four-dimensional subspaces containing $\zeta$  (and, in particular, onto all two-dimensional subspaces containing $\zeta$) coincide or are reflections of each other (with respect to the line containing $\zeta$). Using Theorem 3.1.1 from \cite[page 99]{Ga}  we have 
$\tilde{K}=\tilde{L}$ or  $\tilde{K}={\mathcal O}\tilde{L}$. Thus, 
 $K=L+a_L-a_K$ or $K={\mathcal O}L+{\mathcal O}(a_L)-a_K$.

Now assume that the dimension of $J_1 \cap J_2$ is 2. In this case, let $\{ \zeta,v_1,v_2,v_3\}$ be an orthonormal basis of $J_1$, and  $\{\zeta,v_1,v_2',v_3'\}$ be an orthonormal  basis of $J_2$. Define $J_0$ to be the 4-dimensional subspace with basis  $\{\zeta,v_1,v_2,v_2'\}$.  Then, both $J_1 \cap J_0$ and $J_2 \cap J_0$ have dimension 3, and the above argument can be used. A similar argument can be used if the dimension $J_1 \cap J_2$ is 1. 

Finally, assume that,  in addition, the ``ground" projections $K|G$, $L|G$ onto  all $3$-dimensional subspaces $G$ of $\zeta^{\perp}$, are  
directly congruent and  have no rigid motion symmetries. 
Then, using 
 Theorem \ref{tpr}, we see that $\tilde{K}|J=\tilde{L}|J$ for an arbitrary $4$-dimensional subspace $J$. Hence, the projections of $\tilde{K}$ and $\tilde{L}$
onto all two-dimensional subspaces containing $\zeta$ coincide. Using Theorem 3.1.1 from \cite{Ga} we have
 $\tilde{K}=\tilde{L}$. Thus, $K+a_K=L+a_L$ and
the Corollary is proved. 
\qed

\section{Proofs of Theorem \ref{ts} and Corollary  \ref{nts1}}
\label{sec5}

The proofs are   slightly different from the ones about projections.  
We recall that we consider star-shaped bodies  with respect to the origin.
The direction  $\zeta\in S^3$ will be fixed throughout the proof. 
By the conditions of  Theorem \ref{ts},  the sections $K\cap w^{\perp}$ and $L\cap w^{\perp}$ are directly congruent for every $w\in S^2(\zeta)$.
Hence, for every $w\in S^2(\zeta)$ there exists $\chi_w\in SO(3,S^2(w))$ and $a_w\in w^{\perp}$ such that
\begin{equation}\label{nax211}
\chi_w(K\cap w^{\perp})=(L\cap w^{\perp})+a_w.
\end{equation}

Let $l(\zeta)$ denote the one-dimensional subspace containing $\zeta$. As in Section \ref{sec3}, we use the notation ${\mathcal A}_{K}\subset S^3$  for  the set of directions that are parallel to the diameters of $K$ (similarly for $L$).
We consider the set $\Omega^r$, which is defined similarly to $\Omega$ (see \eqref{nax3}). 
We will use the notation  $v_K(\zeta)=\rho_K(\zeta)+\rho_K(-\zeta)$ for the length of the diameter $d_K(\zeta)$, which contains the origin and is parallel to $\zeta$.

\subsection{Auxiliary lemmata} 
Our first goal is  to reduce matters to rotations leaving $l(\zeta)$ fixed.
We will do this by showing that for  most of the directions $w\in S^2(\zeta)$,
the sections $K\cap w^{\perp}$ and $L\cap w^{\perp}$ have exactly one diameter contained in  $l(\zeta)$.

We will use the well-known properties of the radial function  (see, for example,  \cite[(0.33), page 20]{Ga})
\begin{equation}
 \label{rho}
\rho_{\tilde{K}\cap w^{\perp}}(\theta)=\rho_{\tilde{K}}(\theta),     \; \rho_{\chi_w(\tilde{K}\cap w^{\perp})}(\theta)=\rho_{\tilde{K}\cap w^{\perp}}(\chi_w^{-1}(\theta)), \;  \forall \theta\in w^{\perp}\cap S^3.
\end{equation}

\bl\label{radnu}
Let  $K$ and $L$ be  as in Theorem \ref{ts}.  Then $L$ has a diameter $ d_L(\zeta)$  passing through the origin,  and $\Omega^r$ is everywhere dense in $S^2(\zeta)$. Moreover, for every $w\in \Omega^r$ we have $\chi_w(\zeta)=\pm\zeta$  and $v_K(\zeta)=v_L(\zeta)$.
\el
\bp
Arguing as in the proof of Lemma \ref{nunu2}   (with $\Omega^r$ instead of $\Omega$), we obtain that $\Omega^r$ is everywhere dense in $S^2(\zeta)$,  and that  $\zeta\in {\mathcal A}_L$. 

We will show that there exists a diameter $d_L(\zeta)$ passing through the origin. Assume that this is not true. Then, for each diameter $d_L$ parallel to $\zeta$, the linear subspace $\textrm{span}(d_L)$ is two dimensional. Let ${\mathcal R}(\zeta)$ be the union of all such two-dimensional subspaces, which is a countable union by  the conditions of Theorem \ref{ts}.  Since ${\mathcal A}_L$ is also countable, 
there exists $w\in S^2(\zeta)$ such that $w^{\perp}\cap {\mathcal R}(\zeta)=l(\zeta)$ and  $w^{\perp}$ does not contain any direction $\eta\neq\zeta$ that is parallel to a diameter of $L$.  But then $L$ does not have a diameter in $w^{\perp}$, while $K$ does. This contradiction shows that the diameter in the direction $\zeta$ passes through the origin, and therefore, arguing as in the proof of Lemma \ref{nunu2}, we obtain that 
for all $w\in \Omega^r$, $\chi_w(\zeta)= \pm\zeta$, and $v_K(\zeta)=v_L(\zeta)$. The result follows.
\ep

Our next step is to translate the body $L$ so that its diameter in the direction $\zeta$ coincides with $d_K(\zeta)$. However, if we translate a star-shaped body with respect to the origin, the resulting body need not necessarily be star-shaped with respect to the origin. Lemma \ref{stars} will show that, under the hypothesis of Theorem \ref{ts}, the translation of $L$ must also be star-shaped with respect to the origin (if the bodies $K,L$ are convex, this lemma can be dispensed with). 

\bl\label{stars}
There exists a vector $a \in \mathbb{R}^4$, parallel to $\zeta$, such that the body  $\tilde{L}=L+a$ is star shaped with respect to the origin, and $d_K(\zeta)=d_{\tilde{L}}(\zeta)$. 
\el
\bp
 Consider the sets
\begin{eqnarray*}
    R_1=\{ w \in \Omega^r: \; \chi_w(d_K(\zeta))=d_K(\zeta)\}, \\
    R_2=\{ w \in \Omega^r :\ \;  \chi_w(d_K(\zeta))\neq d_K(\zeta)\},
\end{eqnarray*}
where $\chi_w$ is the rotation in \eqref{nax211}.  
Recall that for any $w\in \Omega^r$, the sections $K \cap w^\perp$ and $L\cap w^\perp$ contain only diameters in the direction $\zeta$. 
If $w \in R_1$, $\chi_w$ must be either a rotation about $\zeta$, or a rotation by angle $\pi$ about some $u \in S^2(\zeta)\cap S^2(w)$ (in this last case,  $d_K(\zeta)$ must be centered at the origin). On the other hand, if $w \in R_2$, $\chi_w$ is a rotation by angle $\pi$ about some $u \in S^2(\zeta)\cap S^2(w)$, and $d_K(\zeta)$ cannot be centered at the origin.

Assume, at first, that $\Omega^r = R_1$. Since the diameter $d_K(\zeta)$ is fixed by  $\chi_w$, the vector $a_w$ in \eqref{nax211} is independent of $w\in \Omega^r$ and  $a_w=a_1=\left(\rho_K(\zeta)-\rho_L(\zeta) \right)\zeta$.
The translated section $(L\cap w^\perp )+a_1$ coincides with $\chi_w(K\cap w^\perp)$, and therefore $(L\cap w^\perp) +a_1$ is star-shaped with respect to the origin for every $w \in \Omega^r$. Since  $\Omega^r$ is dense in $S^2(\zeta)$, we conclude that the translated body $\tilde{L}=L+a$, with $a=a_1$, is also star-shaped with respect to the origin. 

Secondly, assume that  $\Omega^r = R_2$. Then,  $a_w$ is independent of $w \in \Omega^r$ and $a_w=a_2=
\left( \rho_K(-\zeta)-\rho_L(\zeta) \right) \zeta$. We conclude that $\tilde{L}=L+a$, with $a=a_2$, is star-shaped with respect to the origin.

Finally, we show that the case where $R_1$ and $R_2$ are both nonempty does not occur under the assumptions of Theorem \ref{ts}. 
Since $R_1 \cup R_2 =\Omega^r$, we have $S^2(\zeta)=\overline{R_1 \cup R_2} \subseteq \overline{R_1}\cup \overline{R_2} \subseteq S^2(\zeta)$.  Hence, there exists $w_0 \in \overline{R_1}\cap \overline{R_2}$, {\it i.e.}, there is a rotation $\chi_{w_0}$ such that $\chi_{w_0}(d_K(\zeta))=d_K(\zeta)$ and
\begin{equation}\label{r1}
   \chi_{w_0}(K\cap w_0^\perp)=L\cap w_0^\perp +a_1,
\end{equation}
and a rotation $\tilde{\chi}_{w_0}$ such that  $\tilde{\chi}_{w_0}(d_K(\zeta))\neq d_K(\zeta)$ and 
\begin{equation}\label{r2}
   \tilde{\chi}_{w_0}(K\cap w_0^\perp)=L\cap w_0^\perp +a_2.
\end{equation}
 In particular, since $\tilde{\chi}_{w_0}$ does not fix $d_K(\zeta)$, this diameter cannot be centered at the origin, and 
 it follows that  the other rotation $\chi_{w_0}$ must be about $\zeta$. By \eqref{r1} and \eqref{r2} we have
\[
    K\cap w_0^\perp=\chi_{w_0}^{-1} \left( \tilde{\chi}_{w_0}(K\cap w_0^\perp)  \right) + b.
\] 
Observe that the rotation $\chi_{w_0}^{-1} \circ \tilde{\chi}_{w_0}$  is about the vector $v \in S^2(u) \cap S^2(\zeta)$. Since $\chi_{w_0}^{-1} \circ \tilde{\chi}_{w_0}(\zeta)=-\zeta$, this rotation is by angle $\pi$. Therefore, $K\cap w_0^\perp$ has a $(v,\pi)$-symmetry. This contradicts the hypothesis of Theorem \ref{ts}. The Lemma is proven. 
\ep


\bl\label{naxuxu1}
For every $w\in \Omega^r$ there exists $\varphi_w=\chi^{-1}_w\in SO(3, S^2(w))$, $\varphi_w(\zeta)=\pm\zeta$, such that
\begin{equation}\label{naxuxu2}
\rho_{K}\circ \varphi_w(\theta)=\rho_{\tilde{L}}(\theta)\qquad \forall \theta \in S^2(w).
\end{equation}
\el
\bp
In terms of $K$ and $\tilde{L}$, equation \eqref{nax211} can be written as 
\begin{equation}\label{nax3331}
\chi_w(K \cap w^{\perp})=\tilde{L}\cap w^{\perp}
\end{equation}
for some $\chi_w\in SO(3, S^2(w))$, $\chi_w(\zeta)=\pm\zeta$.
Then,
$\rho_{\chi_w(K\cap w^{\perp})}(x)=\rho_{\tilde{L}\cap w^{\perp}}(x)$ for all $x\in w^{\perp}$. In particular, we have that $\rho_{\chi_w(K\cap w^{\perp})}(\theta)=\rho_{\tilde{L}\cap w^{\perp}}(\theta)$ for all $\theta\in S^2(w)$. We now use  (\ref{rho}) to conclude the proof.
\ep

\subsection{Proof of Theorem \ref{ts}} 

Consider the sets
$$
\Xi^r=\{w\in S^2(\zeta):\quad  (\ref{naxuxu2}) \quad\textrm{holds with}\quad \varphi_w(\zeta)=\zeta\}
$$
and 
$$
\Psi^r=\{w\in S^2(\zeta):\quad (\ref{naxuxu2}) \quad\textrm{holds with}\quad \varphi_w(\zeta)=-\zeta\}.
$$

By definition, $\Omega^r \subset (\Xi^r\cup \Psi^r)$. Therefore, Lemma  \ref{radnu} implies that  $\Xi^r\cup \Psi^r=S^2(\zeta)$. 
Now we can apply Proposition \ref{krh2} (with $f=\rho_{K}$, $g=\rho_{\tilde{L}}$, and $\Xi=\Xi^r$, $\Psi=\Psi^r$) 
obtaining that either $\rho_{K}=\rho_{\tilde{L}}$ on  $S^3$, or $\rho_{K}(\theta)=\rho_{\tilde{L}}({\mathcal U}\theta)$  for all $\theta\in S^3$. Here ${\mathcal U}\in O(4)$ is an orthogonal transformation, satisfying ${\mathcal U}|_{S^2(\zeta)}=-I$ and 
${\mathcal U}(\zeta)=\zeta$.
In the first case, $K=\tilde{L}$, and in the second, $K={\mathcal O}\tilde{L}$, where ${\mathcal O}={\mathcal U}^{-1}$. Thus, either $K=L+a$, or $K={\mathcal O}L+b$ with $b={\mathcal O}(a)$.
This finishes the proof of Theorem \ref{ts}.
\qed

\subsection{Proof of Corollary \ref{nts1}}

The proof is similar to the one of Corollary \ref{npr1}. One has only to consider the sections $K \cap J$, $\tilde{L} \cap J$, instead of the projections $K|J$, $\tilde{L}|J$, and Theorem 7.1.1 from  \cite[page 270]{Ga}, instead of Theorem 3.1.1 from  \cite[page 99]{Ga}. 
\qed

\section{Appendix}\label{sec6}
  Let $\delta(K,P)$ be the Hausdorff distance between the convex bodies $K$ and $P$ in ${\mathbb R^n}$, $n\ge 2$,
$\delta(K,P)=\max\limits_{\theta\in S^{n-1}}|h_K(\theta)-h_{P}(\theta)|$. 

Our goal is to prove

\bprop\label{AppProp} Any convex body $K$ in ${\mathbb R^n}$, $n\ge 4$, can be approximated in the Hausdorff metric by polytopes without 3-dimensional projections that have rigid motion symmetries.
\eprop

Since polytopes have finitely many diameters, Proposition \ref{AppProp} shows that the set of bodies satisfying the conditions of Corollary  \ref{npr1}  contains a set of polytopes which is dense in the set of all convex bodies. 

Proposition \ref{AppProp} is not a new  result (see \cite[page 48]{Go1}). An abstract geometric proof of this fact can be given \cite{Pa}. However, for the convenience of the reader, we include an elementary proof.   
The idea is, assuming that $K$ has positive Gaussian curvature, to observe first  that $K$ can be approximated by polytopes whose 3-dimensional projections have many vertices.  If a polytope has a 3-dimensional projection with a rigid motion symmetry, then we use  (\ref{ru1}) to form a system of linear equations, and use the implicit function theorem to prove that these polytopes form a ``manifold" of small dimension.

\subsection{Auxiliary results}
We will need the following theorem and two  lemmata.  Let $C^2_+({\mathbb R^n})$ be the set of convex bodies in ${\mathbb R^n}$ having a positive Gaussian curvature. It is well-known,   that any convex body   can be approximated  in the  Hausdorff metric by  convex bodies $K\in C^2_+({\mathbb R^n})$ \cite[pages 158-160]{Sch}. 
Hence, we can  assume that $K\in C^2_+({\mathbb R^n})$.

Our first auxiliary statement is the following result of  Schneider, \cite{Sch2}.
\bt\label{Mark1}
Let $K\in C^2_+({\mathbb R}^n)$, $n\ge 3$. Then, for $v\to\infty$, we have
$$
\delta(K, P_v^*)\approx c_n\,\,v^{-\frac{2}{n-1}}\,\Big(\int\limits_{\partial K} \sqrt{G_K(\sigma)}d\sigma  \Big)^{\frac{2}{n-1}},
$$
where $P_v^*$ is a polytope with vertices on the boundary  $\partial K$, not unique in general, for which $\delta(K, P_v^*)$  equals the infimum of $\delta(K,P)$ over all convex polytopes $P$ contained in $K$ that have at most $v$ vertices, 
$c_n$ is a constant depending on the dimension, and $G_K(\sigma)$ is the Gaussian curvature of $K$ at $\sigma\in\partial K$.
\et

The next statement is well known. 
\bl\label{Mark3}
Let  $K\in C^2_+({\mathbb R}^n)$, $n\ge 4$.
Then $K|H\in C^2_+(H)$, where $K|H$ is the projection of $K$ onto $H\in {\mathcal G}(n,3)$.
\el
\bp
Let $x$ be any point on the boundary of $K$. Changing  the coordinates if necesssary we can assume that $x$ is the origin and  the tangent hyperplane to $K$ at $x$ is the $(x_1,\dots, x_{n-1})$-hyperplane. Using the  Taylor decomposition of the boundary of $K$ near the origin we have
$$
x_n=f(x_1,\dots,x_{n-1})=k_1x_1^2+\dots+k_{n-1}x_{n-1}^2+o(x),
$$
where $k_j>0$, $j=1,\dots,n-1$, are the main curvatures of the boundary at $x$, and $\frac{o(x)}{|x|}\to 0$ as $|x|\to 0$.
Consider the ball $B$,
$$
B=\{x\in{\mathbb R^n}:\,x_1^2+\dots+x_{n-1}^2+(x_{n}-\frac{1}{k})^2=\frac{1}{k^2}\},\qquad
k=\min\limits_{j=1,\dots,n-1}k_j.
$$
Since the main curvatures are the reciprocals of the main radii of curvature we see that in a small enough neighborhood $W$ of the origin,  $K\cap W$ is contained in $B$. 
Let $u\in S^{n-1}$ be such that  $u_n=0$, {\it i.e.}, $u$ is the unit vector contained in the $(x_1,\dots, x_{n-1})$-hyperplane, and let $H_u\in {\mathcal G}(n,3)$ be 
contained in the $(x_1,\dots, x_{n-1})$-hyperplane, and 
 orthogonal to $u$.
Observe that the boundary of  the projection $(K\cap W)|H_u$
 is contained in the $3$-dimensional ball  of radius $\frac{1}{k}$, which is the projection
of $B$. Since the main curvatures of the boundary of $(K\cap W)|H_u$ are the reciprocals of the radii of curvature, we see that the main curvatures of
$(K\cap W)|H_u$ at the origin are positive. Since $x$ was an arbitrary point on the boundary of $K$, the result follows.
\ep

By Theorem \ref{Mark1}, the number of vertices of $P_v^*$ tends to infinity as $v \rightarrow \infty$.
We claim that for every 3-dimensional space $H$, the number of vertices of the sequence of  polytopes $\{P_v^*|H\}$ is unbounded  as $v \rightarrow \infty$.

\bl\label{vertproj}
Let $V_{H,v}$ be the number of vertices of $P^*_v|H$, and let $\beta_v:=\inf_H V_{H,v}$. Then  the sequence $(\beta_v)_{v=1}^{\infty}$ is unbounded.
\el

\bp
If the sequence $(\beta_v)_{v=1}^{\infty}$ is bounded, then there exists a natural number $m$ such that $\beta_v\le m$ for all $v\in\mathbb{ N}$.
 In particular, for every $v\in \mathbb{R}$ there exists a subspace $H_v\in G(n,3)$ such that
the number of vertices of $P^*_v|H_v$ does not exceed $m$.

Using the first relation in \eqref{nax35}, we have that  $\delta(P^*_v|H_v, K|H_v) \leq \delta(P^*_v, K)$, and that  $ \delta(P^*_v, K)\to 0$ as $v\to\infty$.
For every $v\in {\mathbb N}$, denote by $Q_m=Q_m(v)$ the polytope inscribed in $K|H_v$ such that its distance to $K|H_v$ is minimal among all polytopes inscribed in $K|H_v$ and having at most $m$ vertices.
Then
$\delta(P^*_v|H_v, K|H_v)\ge \delta(Q_m(v), K|H_v)$,
and
$\delta(Q_m(v), K|H_v)\to 0$  as $v\to\infty$.
On the other hand, applying Theorem 3 with $K|H_v$ instead of $K$ and $m$ instead of $v$, we see that
$$
\delta(Q_m(v), K|H_v)\approx m^{-1}\,\int\limits_{\partial (K|H_v)}\sqrt{G_{K|H_v}(\sigma)}d\sigma.
$$
We claim that this is impossible by compactness. Indeed, we can assume that the sequence of bodies $K|H_v$ is convergent, say to $K|H_0$. Hence,
\begin{equation}\label{Maria}
\int\limits_{\partial (K|H_v)}\sqrt{G_{K|H_v}(\sigma)}d\sigma \to
\int\limits_{\partial (K|H_0)}\sqrt{G_{K|H_0}(\theta)}d\theta
\end{equation}
as $v\to \infty$, and
$$
\delta(Q_m(v), K|H_v)\approx m^{-1}\,\int\limits_{\partial (K|H_0)}\sqrt{G_{K|H_0}(\theta)}d\theta.
$$
The left hand side of the last quantity tends to zero, while the right hand side is a positive constant, and we obtain a contradiction.

To show that (\ref{Maria}) holds, we use formula (2.5.29) from   \cite[pg 112]{Sch}, 
$$
\int\limits_{\partial (K|H_v)}\sqrt{G_{K|H_v}(\sigma)}d\sigma=\int\limits_{S^{n-1}\cap H_v}\sqrt{G_{K|H_v}(\nabla h_{K|H_v} (\theta))} {\mathbb{H}}(h_{K|H_v}) (\theta) d\theta.
$$
Here, ${\mathbb H}(h_{K|H_v})$ is the Hessian of the support function (the partial derivatives are the usual derivatives of the support function extended as a  homogeneous function of degree $1$ onto ${\mathbb R}^n$).  
Let $\theta_0\in S^{n-1}\cap H_0$, and consider any sequence $(\theta_v)_{v=1}^{\infty}$, $\theta_v \in S^{n-1}\cap H_v$, converging to $\theta_0$.  Using the fact that $h_{K|H_v}(\theta)=h_K(\theta)$ for $\theta\in S^{n-1}\cap H_v$, we see that as $v\to\infty$ we have
$$
G_{K|H_v}(\nabla h_{K|H_v} (\theta_v))=G_{K|H_v}(\nabla h_{K} (\theta_v))
$$
$$
\to G_{K|H_0}(\nabla h_{K} (\theta_0))=G_{K|H_0}(\nabla h_{K|H_0} (\theta_0)), 
$$
and
$$
{\mathbb H}(h_{K|H_v})(\theta_v)={\mathbb H}(h_{K})(\theta_v)\to {\mathbb H}(h_{K})(\theta_0)={\mathbb H}(h_{K|H_0})(\theta_0).
$$ Hence, (\ref{Maria}) follows.\ep


To formulate our last auxiliary lemma, we recall the definition of the Hausdorff dimension, \cite{WikiH}.
Given any subset $E$ of ${\mathbb R^n}$ and $\alpha\ge 0$,  the {\it exterior $\alpha$-dimensional Hausdorff measure} of $E$ is defined by $m^*_{\alpha}(E)=\lim\limits_{\delta\to 0^+}\inf {\mathcal H}^{\delta}_{\alpha}(E)$, where
$$
{\mathcal H}^{\delta}_{\alpha}(E):=\inf\{\,\,\sum\limits_{k=1}^{\infty}(\,\textrm{diam}\,F_k)^{\alpha}:\,\,E\subset \bigcup\limits_{k=1}^{\infty}F_k,\quad diam\,F_k\le \delta\},
$$
and $diam(S)=\sup\limits_{x,y\in S}|x-y|$ stands for the diameter of $S$. The {\it  Hausdorff dimension} of $E$ is  $dim_H(E)=\inf\{\alpha>0:\,m^*_{\alpha}(E)=0\}$.

\bl\label{Mark2}
Let ${\mathcal M}$ be a smooth manifold of  dimension $k$ in ${\mathbb R}^m$, $m\ge 3$, $k\le m-2$, and let  ${\mathcal M}|H$  be the orthogonal projection of
${\mathcal M}$
onto a $l$-dimensional subspace $H$, $k<l\le m-1$. Then the Hausdorff dimension of ${\mathcal M}|H$ does not exceed the dimension of  ${\mathcal M}$.
\el
\bp
Let $\delta>0$ and let $\bigcup\limits_{k=1}^{\infty}F_k$, $\textrm{diam}(F_k)\le \delta$, be a covering of ${\mathcal M}$. Since 
$\bigcup\limits_{k=1}^{\infty}(F_k|H)$ is a covering of ${\mathcal M}|H$, and 
$\textrm{diam}(F_k|H)\le \textrm{diam}(F_k)\le \delta$, we see that
$$
\sum\limits_{k=1}^{\infty}(\,\textrm{diam}(F_k|H))^{\alpha}\le \sum\limits_{k=1}^{\infty}(\,\textrm{diam}(F_k))^{\alpha},
$$
and $m^*_{\alpha}({\mathcal M}|H)\le m^*_{\alpha}({\mathcal M})$. The result follows.
\ep


\bigskip

\subsection{Proof of Proposition \ref{AppProp}}

To prove the proposition it is enough to show that each $P_v^*$, having sufficiently many vertices, can be approximated by polytopes without any 3-dimensional projection rigid motion symmetries.  We will do this by proving that the set of polytopes having $v$ vertices {\it with} 3-dimensional projection rigid motion symmetries is a nowhere dense set contained in the set of all polytopes having $v$ vertices.

Define ${\mathcal P}_v$ to be 
the set of polytopes in ${\mathbb R^n}$, $n\ge 4$, with $v$ vertices  $p_1, p_2,\dots, p_v$.
We see that ${\mathcal P}_v$ can be parametrized by points from ${\mathbb R^{nv}}$, with
$p_j=(p_{1j},\dots, p_{nj})\in {\mathbb R^n}$, $j=1,\dots, v$, and  we can identify ${\mathcal P}_v$ with an open domain in ${\mathbb R}^{nv}$.

We denote by  $\Pi_v$   the set of  polytopes in ${\mathcal P}_v$ that have a  3-dimensional projection with  rigid motion symmetries. 
Our goal is to show that $\Pi_v$  is  nowhere dense in ${\mathcal P}_v$, provided that $v$ is large enough.  We can partition $\Pi_v$ into equivalence classes such that two polytopes are in the same class if there is a rigid motion in ${\mathbb R}^n$ taking one to the other.  Letting $H_0$ be the $(x_1,x_2,x_3)$-plane in ${\mathbb R}^n$, each equivalence class can be represented by a polytope whose projection on $H_0$ has rigid motion symmetries.  Let us define ${\mathcal Q}_v$ to be the set of these representatives, {\it i.e.}, 
\[ {\mathcal Q}_v  
= \{ Q \in {\mathcal P}_v : \exists \varphi_{H_0} \in O(3,H_0), \varphi_{H_0} \not= I, \exists a_{H_0} \in {\mathbb R}^3 \textrm{ such that} 
\]
\begin{equation}\label{ru1}
\qquad \varphi_{H_0}(Q|H_0)+a_{H_0}=Q|H_0 \}.
\end{equation}

Observe that every $P \in \Pi_v$ can be written as $P=\phi(Q)+b$ for some $\phi \in O(n)$, $Q \in {\mathcal Q}_v$, $b \in {\mathbb R}^n$, and hence can be represented as the triple $(Q,\phi,b) \in {\mathcal Q}_v \times O(n) \times {\mathbb R}^n$.  Thus, 
\begin{equation}\label{dim}
\dim (\Pi_v) \leq \dim ({\mathcal Q}_v) + \dim (O(n)) + n = \dim ({\mathcal Q}_v) + \frac{n(n+1)}{2}.
\end{equation}
All that remains is to find the dimension of ${\mathcal Q}_v$.  
Consider the set ${\mathcal M}={\mathcal M}({\mathcal Q}_v)$ of all triples 
$$
(Q, \varphi_{H_0}, a_{H_0}) \in {\mathbb R^{nv}} \times O(3,H_0) \times {\mathbb R^3},
$$
satisfying  (\ref{ru1}).

Let $H \in {\mathcal G}(n,3)$. By \eqref{nax35}, for every $\theta\in H\cap S^{n-1}$, we have $h_{K|H}(\theta)=h_K(\theta)$. Thus,  $K|H$ can be approximated in the Hausdorff metric by polytopes $P_v^*|H$.  By Lemma  \ref{vertproj}, we have a subsequence $\{\beta_{v_j}\}$ of $\{\beta_v\}$, such that  
 \begin{equation}
   \label{vertices}
 \beta_{v_j}>5+\frac{n(n+1)}{2}
  \end{equation}
  if $j$ is large enough. For simplicity, in the following we will denote $\beta_{v_j}$ by $\beta_v$.
\bl\label{ver1}
The set ${\mathcal M}$ is a 
manifold in ${\mathbb R}^{nv+6}$ with dimension at most $(nv+5-\beta_v)$, provided that 
$v$ is such that $\beta_v>5+\frac{n(n+1)}{2}$. 
\el
\bp
Let $Q$ be a polytope in ${\mathcal Q}_v$ and consider its projection $Q|H_0$, which is also a polytope with $t$ vertices, where $t \geq \beta_v$. We will  write the assumption that $Q|H_0$ has rigid motion symmetries as a system of  linear equations that equal zero precisely at the vertices of $Q|H_0$, and explicitly compute the determinant of its Jacobian matrix to show that it is nonzero. The implicit function theorem \cite{Wiki} will allow us to obtain the result. 

Since any rigid motion maps a vertex into a vertex, an equation, similar to (\ref{ru1}), can be written for the corresponding vertices $q_i|H_0$ of $Q|H_0$,
\begin{equation}\label{ru2}
q_i|H_0=\varphi_{H_0}(q_{j(i)}|H_0)+a_{H_0},
\end{equation}
where  $\varphi_{H_0}$ is a nonidentical orthogonal transformation whose $3 \times 3$ matrix has coordinates  $(o_{l,m})_{l,m=1,2,3}$, 
and $j$ is a permutation on the set  $\{1,\ldots t\}$, which indicates that  the $j(i)$-th vertex gets mapped to the $i$-th vertex.  As it is well known, a permutation can be written as a product of cycles. We will consider two cases: cycles of length one, and cycles of length greater than one.  

Assume that the vertex $q_i|H_0$ is mapped to itself, {\it i.e.},  $q_i|H_0=\varphi_{H_0}(q_i|H_0)+a_{H_0}$.
   Since $\varphi_{H_0}$ is not the identity,  given a basis $e_1,e_2,e_3$ of $H_0$, there exists $r \in \{1,2,3\}$ such that $\varphi_{H_0}(e_r) \not= e_r$.  For this $r$, consider the function 
$F_{ri}:\mathbb{R}^{nv} \times O(3,H_0) \times \mathbb{R}^3 \rightarrow \mathbb{R}$ defined by 

$$F_{ri}(x_{11},\ldots,x_{nv},\varphi_{H_0},a_{H_0})=((x_{1i},x_{2i},x_{3i})-\varphi_{H_0}(x_{1i},x_{2i},x_{3i})-a_{H_0})_r.$$
\[
   =x_{ri}-o_{r1}x_{1i}-o_{r2}x_{2i}-o_{r3}x_{3i}-(a_{H_0})_r.
\]
Since the right hand side depends only on the variables $x_{1i},x_{2i},x_{3i}$, we see that 
$\frac{\partial F_{ri}}{\partial x_{ks}}=0$ for all $s \not= i$ and all $k$, while $\frac{\partial F_{ri}}{\partial x_{ri}} \not= 0$ because $\varphi_{H_0}(e_r) \not= e_r$. Thus, this cycle forms a $(1\times 1)$-Jacobian block whose entry is not 0.

Next, suppose that the cycle is of length $k$ and permutates the vertices
$q_{i_1},q_{i_2},\ldots,q_{i_k}$, (for $\ell<k$, $q_{i_{\ell+1}}$ gets mapped to $q_{i_{\ell}}$ and $q_{i_1}$ is mapped back to $q_{i_k}$).  Consider the system of $3(k-1)$ functions $F_{rs}:\mathbb{R}^{nv} \times O(3,H_0) \times \mathbb{R}^3 \rightarrow \mathbb{R}$ defined by 
$$
F_{rs}(x_{11},\ldots,x_{nv},\varphi_{H_0},a_{H_0})=((x_{1s},x_{2s},x_{3s})-\varphi_{H_0}(x_{1j(s)},x_{2j(s)},x_{3j(s)})-a_{H_0})_r
$$
for  $r=1,2,3$ and for  $s=i_1,i_2,\ldots,i_{k-1}$.

We will order the variables in such a  way that the Jacobian block corresponding to this cycle will be upper  triangular.  We note that for $r=1,2,3$, and $s=i_1,\ldots,i_{k-1}$, $F_{rs}$ depends on the variables $x_{rs}$ and $x_{kj(s)}$ for $k=1,2,3$.  Thus, $\frac{\partial F_{rs}}{\partial x_{ks}}=0$ for $k \not= r$, and $\frac{\partial F_{rs}}{\partial x_{k\ell}}=0$ for all $\ell \not= s$, $\ell \not= j(s)$ and all $k$.  Order the Jacobian block as follows, $x_{1i_1},x_{2i_1},x_{3i_1},x_{1i_2},\ldots,x_{3i_{k-1}}$.  Since $\frac{\partial F_{rs}}{\partial x_{rs}}=1$, the diagonal entries are all 1.  In addition, the variables $x_{kj(s)}$ occur after $x_{rs}$, so the Jacobian block is upper triangular.  Therefore, the determinant of this block is equal to 1. Thus,  the Jacobian of the system of equations  is a block diagonal matrix with nonzero determinant.




We observe that the number of equations in our system depends on the decomposition of the permutation $j$ into cycles. Each 1-cycle gives us one equation, while each cycle of length $k>1$ contributes $3(k-1)$ equations to the system.  Hence, the smallest possible number of equations in our system is $3+(t-2)$, which occurs if the decomposition of the permutation  $j$ into cycles contains only one two-cycle and all the rest are one-cycles.  By the implicit function theorem, we can express at least $t+1$ variables $x_{rs}$ as functions of the coordinates of $\varphi_{H_0}, a_{H_0}$  and at most $nv-(t+1)$ other variables. Since $t \geq \beta_v$, this shows that the dimension of the manifold  ${\mathcal M}$  in ${\mathbb R}^{nv+6}$.  is at most  $\left( nv+\dim(O(3)+\dim(H_0)-(\beta_v+1) \right)=nv+5-\beta_v$.
\ep

We are now ready to prove our goal.

\bl\label{Mark4}
The set $\Pi_v$ is  nowhere dense in ${\mathcal P}_v$.
\el
\bp
By definition,  ${\mathcal Q}_v$ is equal to the projection of ${\mathcal M}$ onto ${\mathbb R}^{nv}$ 
and by Lemmata \ref{Mark2} and  \ref{ver1}, 
$$
\dim({\mathcal Q}_v)=\dim({\mathcal M}|{{\mathbb R}^{nv}}) \leq \dim({\mathcal M}) \leq nv+5-\beta_v.
$$
Hence, using  (\ref{dim}), we have $\dim(\Pi_v) \leq nv+5-\beta_v+\frac{n(n+1)}{2}$.  Finally, (\ref{vertices}) yields 
 $\dim(\Pi_v)<\dim({\mathcal P}_v)=nv$.
\ep

To complete the proof of  Proposition \ref{AppProp}, we use  Theorem \ref{Mark1} to approximate $K\in C^2_+({\mathbb R^n})$ in the Hausdorff metric, by polytopes $P_v^*$ with $v$ so large that $t_0>5+\frac{n(n+1)}{2}$.  By Lemma \ref{Mark4}, we can approximate $P_v^*$ by polytopes without 3-dimensional projections that have rigid motion symmetries. \hfill $\square$ \\

\end{document}